\input ./style/arxiv-general.cfg
\documentclass[MSNbibl,number,seceqn,citesort,dvips]{arxbj}
\makeatletter
   \@ifpackageloaded{graphicx}{}{\usepackage{graphicx}}
\makeatother

\volume{22}
\issue{1}
\pubyear{2016}
\firstpage{530}
\lastpage{562}
\doi{10.3150/14-BEJ667} 
\docsubty{FLA}

\makeatletter
\def\sfrac#1#2{#1/#2}
\def\vfrac#1#2{(#1)/#2}
\def\afrac#1#2{#1/(#2)}

\newcommand{\rrvert}{\vert}

\newcommand{\rrVert}{\Vert}
\newcommand{\llvert}{\vert}
\newcommand{\llVert}{\Vert}
\newproclaim{definition}{Definition}[section]
\newtheorem{theorem}[definition]{Theorem}
\newtheorem{proposition}[definition]{Proposition}
\newtheorem{lemma}[definition]{Lemma}
\newtheorem{corollary}[definition]{Corollary}
\DeclareMathAlphabet{\mathpzc}{OT1}{pzc}{m}{it}
\newcommand{\E}{\mathbb{E}}
\newcommand{\var}{\mathbb{V}\mathrm{ar}}
\renewcommand{\L}{\mathbf{L}}
\renewcommand{\P}{\mathbb{P}}
\newcommand{\R}{\mathbb{R}}
\newcommand{\N}{\mathbb{N}}
\newcommand{\vecx}[1]{\mathbf{\underline{#1} } } 
\newcommand{\cX}{\mathcal{X}}
\newcommand{\1}{\mathbf{1}}
\newcommand{\thetaL}{\theta_L}
\newcommand{\thetaC}{\theta_c}
\newcommand{\thetaP}{\theta_\Phi}
\newcommand{\thetaconv}{\theta_{\mathrm{conv}}}
\newcommand{\intB}[2]{B_{#1,#2}}
\newcommand{\Cyi}{C_{y,i}}
\newcommand{\Cyj}{C_{y,j}}
\newcommand{\Cyjp}{C_{y,j+1}}
\newcommand{\Cyk}{C_{y,k}}
\newcommand{\Czi}{C_{z,i}}
\newcommand{\Czj}{C_{z,j}}
\newcommand{\aprioriConst}[2]{C^{(#1)}_{#2} }
\newcommand{\cyz}{C_{y,z}}
\newcommand{\cw}{{\mathcal C}_{(\ref{eqiterationfeed2}\mathrm{a})}}
\newcommand{\cuh}{{\mathcal C}_{(\ref{eqiterationfeed2}\mathrm{b})}}
\newcommand{\cwg}[1]{{\mathcal C}^{(#1)}_{(\ref{eqintuw})}}
\newcommand{\cuv}{{\mathcal C}_{\mathrm{(\ref{equpperboundu})}}}
\newcommand{\yM}[2]{y^{(M)}_{#1}(#2)}
\newcommand{\yMa}[2]{y^{(M)}_{#1}\bigl(#2\bigr)}
\newcommand{\zM}[2]{z^{(M)}_{#1}(#2)}
\newcommand{\zMa}[2]{z^{(M)}_{#1}\bigl(#2\bigr)}
\newcommand{\XM}[2]{X^{(#1,#2)}} 
\newcommand{\MWmarkov}[2]{ \mathpzc{h}^{(#1)}_{#2} }
\newcommand{\Xmarkov}[2]{ V^{(#1)}_{#2} }

\newcommand{\pithetapi}{\pi^{(\theta_\pi)}}
\newcommand{\pit}[1]{\pi^{(#1)}}
\newcommand{\tpi}{\theta_\pi}

\newcommand{\X}[2]{X^{(#1)}_{#2}}
\newcommand{\MW}[2]{H^{(#1)}_{#2}}
\newcommand{\Samp}[1]{\cX^{(#1)} }
\newcommand{\ObsM}[2]{S^{(M)}_{#1}(#2)}
\newcommand{\ObsMa}[2]{S^{(M)}_{#1}\bigl(#2\bigr)}
\newcommand{\Obs}[2]{S_{#1}(#2)}
\newcommand{\Obsa}[2]{S_{#1}\bigl(#2\bigr)}
\newcommand{\tObsM}[2]{\tilde S^{(M)}_{#1}(#2)}
\newcommand{\tObsMa}[2]{\tilde S^{(M)}_{#1}\bigl(#2\bigr)}
\newcommand{\xiM}[2]{\xi^*_{#1}(#2)}
\newcommand{\txiM}[2]{\tilde\xi^*_{#1}(#2)}
\newcommand{\ObsBd}[1]{\bar C_{#1} }

\newcommand{\plx}[3]{p^{(#2)}_{#1}(#3) }

\newcommand{\pl}[2]{p_{#1}^{(#2)} }
\newcommand{\OLS}{\mathrm{OLS}}

\newcommand{\locpol}[1]{{\mathcal P}_{\mathrm{loc.}}^{#1}}
\newcommand{\F}[2]{\mathcal{F}^{(#1)}_{#2}}
\newcommand{\Om}[1]{\Omega^{(#1)} }
\newcommand{\Prob}[1]{\P^{(#1)} }

\newcommand{\HHp}{{$(\mathbf{A'_H})$}}
\newcommand{\HF}{{$(\mathbf{A_F})$}}
\newcommand{\HFc}{{$(\mathbf{A_F})$\textup{(iii)}}}
\newcommand{\HFa}{{$(\mathbf{A_F})$\textup{(i)}}}
\newcommand{\HFb}{{$(\mathbf{A_F})$\textup{(ii)}}}
\newcommand{\HFab}{{$(\mathbf{A_F})$\textup{(i)--(ii)}}}
\newcommand{\HG}{{$(\mathbf{A}_{\bolds{\xi}})$}}
\newcommand{\HH}{{$(\mathbf{A_H})$}}
\newcommand{\HX}{{$(\mathbf{A_X})$}}
\newcommand{\HFp}{{$(\mathbf{A'_F})$}}
\newcommand{\HGp}{{$(\mathbf{A'_{\bolds{\xi}}})$}}
\newcommand{\HGpa}{{$(\mathbf{A'_{\bolds{\xi}}})$\textup{(i)}}}
\newcommand{\HGt}{{$(\mathbf{A''_{\bolds{\xi}}})$}} 
\makeatother

\begin{document}
\begin{frontmatter}

\title{Approximation of backward stochastic differential equations using Malliavin weights and least-squares regression}
\runtitle{Approximation of BSDEs using Malliavin weights and
least-squares regression}

\begin{aug}
\author[A]{\inits{E.}\fnms{Emmanuel}~\snm{Gobet}\thanksref{e1}\ead
[label=e1,mark]{emmanuel.gobet@polytechnique.edu}}
\and
\author[A]{\inits{P.}\fnms{Plamen}~\snm{Turkedjiev}\corref{}\thanksref
{e2}\ead[label=e2,mark]{turkedjiev@cmap.polytechnique.fr}}
\address[A]{Centre de Math\'ematiques Appliqu\'ees,
Ecole Polytechnique and CNRS,
Route de Saclay,
91128 Palaiseau cedex, France.\\
\printead{e1,e2}}
\end{aug}

\received{\smonth{8} \syear{2013}}
\revised{\smonth{3} \syear{2014}}

%
\begin{abstract}
We design a numerical scheme for solving a Dynamic Programming equation
with Malliavin weights arising from
the time-discretization of backward stochastic differential equations
with the integration by parts-representation
of the $Z$-component by
(\textit{Ann. Appl. Probab.} \textbf{12} (2002) 1390--1418). When the sequence of
conditional expectations is computed using
empirical least-squares regressions, we establish, under general
conditions, tight error bounds as the time-average
of local regression errors only (up to logarithmic factors).
We compute the algorithm complexity by a suitable optimization of the
parameters, depending on the dimension
and the smoothness of value functions, in the limit as the number of
grid times goes to infinity. The estimates
take into account the regularity of the terminal function.
\end{abstract}

%
\begin{keyword}
\kwd{backward stochastic differential equations}
\kwd{dynamic programming equation}
\kwd{empirical regressions}
\kwd{Malliavin calculus}
\kwd{non-asymptotic error estimates}
\end{keyword}
\end{frontmatter}

\section{Introduction}
\subsection{Setting}
Let $T>0$ be a fixed terminal time and let $(\Omega, \mathcal
{F},(\mathcal{F}
_t)_{0\leq
t \leq T},\P)$ be a filtered probability space whose filtration is
augmented with the $\P$-null sets.
Let $\pi= \{ 0 =: t_0 < t_1 < \cdots< t_{N-1} < t_N:=T\}$ be a given
time-grid on $[0,T]$ and $\Delta_i:=t_{i+1}-t_i$.
Additionally, for a fixed $q\in\N\setminus\{0\}$, we are given a set
$\{ \MW ij\dvt  0\le i < j \le N\}$ of $(\R^q)^\top$-valued random
variables in $\L_2(\mathcal{F}_T,\P)$ (i.e., square integrable and
$\mathcal{F}
_T$-measurable) that we call Malliavin weights. Here ${}^\top$ stands
for the transpose.

In this paper, we introduce the Malliavin Weights Least Squares
algorithm, abbreviated \textit{MWLS}, to approximate discrete time
stochastic processes $(Y,Z)$
defined by
%
\begin{eqnarray}
\cases{
\displaystyle Y_i =
\mathbb{E}_i\Biggl[\xi+ \sum_{j=i}^{N-1}
f_j(Y_{j+1},Z_j) \Delta_j
\Biggr], &\quad $0 \le i \le N$,
\vspace*{5pt}\cr
\displaystyle Z_i= \mathbb{E}_i\Biggl[\xi\MW iN +
\sum_{j=i+1}^{N-1} f_j(Y_{j+1},Z_j)
\MW ij \Delta_j\Biggr], &\quad $0 \le i \le N-1$,}\label{eqmalscheme}
\end{eqnarray}
where $\mathbb{E}_i[\cdot]:= \mathbb{E}[\cdot\mid\mathcal{F}_{t_i}]$,
$\xi
$ is a $\R$-valued
random variable in $\L_2(\mathcal{F}_T,\P)$, and $(\omega,y,z)
\mapsto
f_j(\omega,y,z)$ is $\mathcal{F}_{t_j} \otimes\mathcal{B}(\R
)\otimes
\mathcal{B}((\R
^q)^\top
)$-measurable. This system is solved backward in time in the order
$Y_N, Z_{N-1}, Y_{N-1},\dots$ and it takes the form of a dynamic
programming equation with Malliavin weights. We call it the Malliavin
Weights Dynamic Programming equation (MWDP for short).

The main application of (\ref{eqmalscheme}) is to approximate
continuous-time, {decoupled} Forward--Backward SDEs (FBSDEs) of the form
%
\begin{equation}
Y_t = \xi+ \int_t^T
f(s,X_s,Y_s,Z_s) \,\mathrm{d}s - \int
_t^T Z_s \,\mathrm{d}W_s,
\label{eqBSDEcts}
\end{equation}
where $(W_s)_{s\ge0 }$ is a Brownian motion in $\R^q$, $(X_s)_{s\ge0
}$ is a diffusion in $\R^d$ and $\xi$ is of the form $\Phi(X_T)$.
Indeed, the MWDP (\ref{eqmalscheme}) was inspired by \cite{mazhan02}, Theorem~4.2, which states that there is a version of the
{continuous-time} process $(Z_t)_{0\le t<T}$ given by
%
\begin{equation}
Z_t = \mathbb{E}_t\biggl[\xi\MW{t}T + \int
_{t}^T f(s,X_s,Y_s,Z_s)
\MW{t}s \,\mathrm{d}s\biggr], \label{eq4Z}
\end{equation}
where the processes $(\MW ts)_{0\leq t < s \leq T} $ are Malliavin
weights defined by
%
\begin{equation}
\MW ts= \frac{1}{s-t} \biggl( \int_t^s
\bigl(\sigma^{-1}( r,X_r) D_tX_r
\bigr)^\top \,\mathrm{d}W_r \biggr)^\top, \qquad0 \le t< s
\le T, \label{eq4malweights}
\end{equation}
for $(D_tX_r)_t$ being the Malliavin derivative of $X_r$ and $\sigma
(\cdot)$ is the diffusion coefficient of $X$. The representation
(\ref{eq4Z}) is obtained via a Malliavin calculus integration by parts
formula, see \cite{nual95} for a general account on the subject.
A discretization procedure to approximate the FBSDE \mbox{(\ref{eqBSDEcts})--(\ref{eq4Z})} with (\ref{eqmalscheme}),
including explicit definitions of the random variables $\MW ij$ based
on~(\ref{eq4malweights}),
is given in \cite{turk13,turk13b}, where the author also
computes the discretization error in terms of $N$.
In honour of the connection between (\ref{eqmalscheme}) and (\ref
{eqBSDEcts})--(\ref{eq4Z}), we call the random variables $\MW ij$
Malliavin weights, $\xi$ the terminal condition, and $(i,\omega,y,z)
\mapsto f_i(y,z)$ the driver.
We say that the pair $(Y,Z)$ satisfying (\ref{eqmalscheme}) solves a
MWDP with
{terminal condition $\xi$ and driver $f_i(y,z)$.}

\subsection{Contributions}
In this paper, we are not concerned with the discretization procedure,
but rather with the analysis of the MWDP equation (\ref{eqmalscheme})
itself and
its numerical resolution via the MWLS algorithm, in which one uses
empirical least-squares regressions (approximations on finite basis of
functions using simulations) to compute conditional expectations.
Since the system (\ref{eqmalscheme}) may be relevant to problems
beyond the FBSDE system (\ref{eqBSDEcts})--(\ref{eq4Z}), we allow the
framework and assumptions to
accomodate as much generality as possible. However, MWLS is, to the
best of our knowledge, the first direct implementation of formula
(\ref{eq4Z}) in a fully implementable numerical scheme.
For other applications of Malliavin calculus in numerical simulations,
with different perspectives and results to ours, see, for instance, \cite
{fourlion99,kohapett02,bouctouz04,gobemuno05,ballcarazane05,hunualsong11,brialaba13}.

We adapt the recent theoretical analysis of the
Least Squares Multi-step forward Dynamical Programming algorithm (LSMDP) of \cite
{gobeturk14} for discrete BSDEs (without Malliavin weights) to the
setting of MWDP.
As in the aforementioned reference, we consider a locally Lipschitz
driver $f_i(y,z)$ that is locally bounded at $(y,z) = (0,0)$ -- see
Section~\ref{section4assumptions}.
This allows the algorithm of the current paper to be applied for the
approximation of some quadratic BSDEs and for some proxy/variance
reduction methods.
For more details on these applications, see \cite{gobeturk14}, Section~2.2.
Moreover, we make use of analogous stability results and conditioning
arguments in the proof of the main result, Theorem~\ref{thmMCerr}, as
in the proof of \cite{gobeturk14}, Theorem 4.11.
{However, the Malliavin weights lead to significantly differences in
the main theorem and stability results, both in the technical elements
of the proofs and the results. We develop {seemingly} novel Gronwall
type inequalities to handle the technical differences; these results
are outlined in Section~\ref{sectiongronwall} and proved in
Appendix~\ref{appendixprooflemintegrals}.
Furthermore, the stability results are more powerful and the complexity
of the MWLS is better than the LSMDP, as will be discussed in what follows.
}

{We would like to mention that the class of quadratic problems we can
treat with these assumptions is quite different to the recent \cite
{chasrich13}.
Here we are treating the {non-degenerate} setting where the terminal
condition may be H\"older continuous, whereas the other reference
allows degeneracy at the expense of requiring locally Lipschitz
terminal conditions.
}

We prove stability results on the MWDP in Section~\ref{sectionstability}.
{Much effort is made to keep the constants explicit.}
These results are instrumental throughout the paper. The stability
estimates on $Z$ are at the individual time points (coherently with the
representation theorem of \cite{mazhan02})
rather than the time-averaged estimates of \cite{gobeturk14}, Proposition~3.2.
This allows for finer and more precise computations.
The time-dependency in our estimates also takes better into account the
regularity of the terminal condition, similarly to the continuous-time
case \cite{delaguat06}.

Section~\ref{section4MC} is the core of the paper: it is dedicated to
the MWLS algorithm in the Markovian context $Y_i=y_i(X_i)$ and
$Z_i=z_i(X_i)$ for some Markov chain $X_i$ in $\R^d$ and unknown
functions $(y_i(\cdot),z_i(\cdot))$. In MWLS, the conditional
expectations in (\ref{eqmalscheme})
are replaced by Monte Carlo least-squares regressions. For each point
of the time-grid, we use a cloud of independent paths of the {Markov
chain} $X$ and the Malliavin weights $H$, and some approximation spaces
for representing the value functions $(y_i(\cdot),z_i(\cdot))$. The
algorithm is detailed in Section~\ref{sectionMCalg}, and a full
error analysis
in terms of the number of simulations and the approximation spaces
is performed in Sections~\ref{sectionerr} and~\ref
{sectionproofmainthm}. The final error estimates (Theorem~\ref
{thmMCerr}) are similar to \cite{gobeturk14}, Theorem 4.11, in that
they are the time-averaged {local} regression errors of
the discrete BSDE, but the results are in a stronger norm and the
time-dependency is better. {The constants are completely explicit.}
{Although the norms are stronger than in \cite{gobeturk14}, the
estimates do not deteriorate; instead, they are {significantly}
improved. This is intrinsically due to the MWDP representation, which
avoids the usual \mbox{$1/\Delta_i$-}factor in front of all controls on $Z$.
This improvement can be tracked by inspecting the a.s. bounds (compare
(\ref{eqasz}) and \cite{gobeturk14}, equation~(14)) and the statistical
error bounds (compare the $\frac{ K_{Z,k} }{M_k }$-terms in (\ref
{eqerr}) of Theorem~\ref{thmMCerr} and the $\frac{ K_{Z,k}
}{\Delta_k
M_k }$-terms of \cite{gobeturk14}, Theorem 4.11).
These error estimates {are} optimal with respect to the convergence
rates (up to logarithmic factors)
under high generality regarding the distribution of the stochastic
model for $X$ and $H$, even if the constants may be conservative. This
is because the local regression errors are optimal under model-free
estimates (Proposition~\ref{propeqyzerrdecoM}).}

{
With the error estimates of Theorem~\ref{thmMCerr} in hand, we
perform a complexity analysis in Section~\ref{sectioncomplexity}. We
propose a choice of basis functions and use it to calibrate the number
of simulations in order to achieve a specified error level.
This then allows us to compute the complexity of the algorithm for that
error level.
The methodology for doing this is analogous to \cite{gobeturk14}, Section~4.4, in that we use the same basis functions
-- which also enable us to study the benefit of smoothness properties of
the underlying Markov functions $(y_i(\cdot),z_i(\cdot))$ --
and also in that we set the ensure the global error level by
calibrating the local regression errors.
However, the conclusion of this section is that MWLS yields better
performance in terms of complexity than LSMDP.
The main reason for this is the improved time-dependancy of the error
estimates, which is a systemic improvement that allows one to make
generate fewer simulations to obtain certain error levels.
Unfortunately, the complexity reduction does not reduce the dependence
on the dimension compared to the LSMDP.
The curse of dimensionality still appears, and the rates depend on the
dimension of the Markov chain $X$ (i.e., $d$).
Nevertheless, the reduction of complexity is substantial and, since one
is able to make fewer simulations to obtain the same error level, will
help alleviate the pressure on memory resources that is typical with
{least-squares} Monte Carlo algorithms.}

This paper is theoretically oriented, and is aimed at paving the way
for new such numerical approaches {based on Malliavin calculus}. Future
works will be devoted to a deeper investigation about the numerical
performance of the MWLS algorithm compared to other known numerical schemes.

\subsection{Notation used throughout the paper}\label{sectionnotation}
\begin{itemize}
\item$\llvert x\rrvert$ stands for the Euclidean norm of the
vector $x$, {${ }^\top
$ denotes the transpose operator}.
\item$\llvert U\rrvert_{p}:=(\E[\llvert U\rrvert
^p])^{{1}/{p }}$ stands for the $\L
_p(\P
)$-norm ($p\geq1$) of a random variable $U$. The \mbox{$\mathcal{F}
_{t_k}$-}conditional version is denoted by $\llvert U\rrvert
_{p,k}:=(\E
_k[\llvert U\rrvert^p])^{1/p}$.
To indicate that $U$ is additionally measurable w.r.t. the $\sigma
$-algebra $\mathcal{Q}$, we may write $U\in\L_{ p}(\mathcal{Q},\P)$.
\item For a multi-dimensional process $U=(U_i)_{0\leq i\leq N}$, its
$l$th component is denoted by $U_l=(U_{l,i})_{0\leq i\leq N}$.
%
\item For any finite $L>0$ and $x = (x_1,\ldots,x_n)\in\R^n$, define
the truncation function
$\mathcal{T}_L(x):= (-L\vee x_1 \wedge L,\ldots, -L \vee x_n \wedge L)$.
\item For finite $x >0$, $\log(x)$ is the natural logarithm of $x$.
\end{itemize}

\subsection{Assumptions}\label{section4assumptions}
\subsubsection*{First set of hypotheses} The following assumptions hold
throughout the entirety of the paper.
Let $R_\pi>0$ be a fixed parameter: this constant determines which
time-grid can be used. The larger $R_\pi$, the larger the class of
admissible time-grids. All subsequent error estimates depend on $R_\pi$.
\begin{enumerate}[\HF]
\item[\HG] $\xi$ is in $\L_2(\mathcal{F}_T,\P)$.
\item[\HF]
\begin{enumerate}[(iii)]
\item[(i)] $(\omega,y,z)\mapsto f_i(y,z)$ is $\mathcal
{F}_{t_i}\otimes
\mathcal
{B}(\R) \otimes\mathcal{B}((\R^{q})^\top)$-measurable for every
{$i<N$}, and
there exist deterministic parameters $\thetaL\in(0,1]$ and $L_{f}\in
[0,+\infty)$ such that
%
\[
\bigl\llvert f_i(y,z)-f_i
\bigl(y',z'\bigr)\bigr\rrvert\leq\frac{ L_f
}{(T-t_i)^{(1-\thetaL
)/2}}
\bigl(\bigl\llvert y-y'\bigr\rrvert+\bigl\llvert z-z'
\bigr\rrvert\bigr),
\]
for any $(y,y',z,z')\in\R\times\R\times(\R^q)^\top\times(\R
^q)^\top$.
\item[(ii)] There exist deterministic parameters $\thetaC\in(0,1]$ and
$C_{f}\in[0,+\infty)$ such that
%
\[
\bigl\llvert f_i(0,0)\bigr\rrvert\leq
\frac{ C_f }{(T-t_i)^{1-\thetaC
}}, \qquad\forall0 \le i < N.
\]
\item[(iii)] The time-grid $\pi:= \{0=t_0 < \cdots< t_N = T \}$ satisfies
\[
\max_{0\le i\le N-2}\frac{\Delta_{i+1}}{\Delta_{i}} \le R_\pi.
\]
\end{enumerate}
\item[\HH] For all $0\le i <j \le N$, the Malliavin weights satisfy
%
\[
\mathbb{E}\bigl[\MW ij\mid\mathcal{F}_{t_i}\bigr] = 0, \qquad
\bigl(\mathbb{E}\bigl[\bigl\llvert\MW ij\bigr\rrvert^2\mid
\mathcal{F}_{t_i}\bigr] \bigr)^{1/2}\le\frac{C_{M}
}{(t_j-t_i)^{1/2}}
\]
for a finite constant $C_{M}\geq0$.
\end{enumerate}
\textit{Comments.} We remark that assumptions \HG~and \HFab~are the
same as their equivalents in \cite{gobeturk14}, Section~2. The usual
case of ``Lipschitz'' BSDE is covered by $\thetaL=\thetaC=1$. As
explained in \cite{gobeturk14}, the case of locally Lipschitz driver
($\thetaL<1$ and/or $\thetaC<1$) is interesting because it allows a
large variety of applications, such as solving BSDEs using proxy
methods or variance reduction methods, and solving quadratic BSDEs. We
refer the reader to \cite{gobeturk14}, Section~2.2, for details.

Assumption \HFc~is much simpler compared to \cite{gobeturk14}. If
$R_\pi\geq1$, \HFc~is satisfied by any time grid with non-increasing
time-step, such as the grids of \cite
{gobemakh10,rich11,geisgeisgobe12}. This may be valuable for
future work on time-grid optimization.

Assumption \HH~is specific to the dynamic programming equation with
Malliavin weights. It is satisfied for the weights derived in \cite
{mazhan02}, and this can remain true after discretization (see \cite
{turk13} or \cite{gobemuno05}).
{
It is also satisfied if $X$ takes the form $X_t=g(t,W_t)$ (like
multi-dimensional geometric Brownian motion), by a simple change of
variables one can use the Malliavin weights $\MW ts= \frac
{(W_s-W_t)^\top}{s-t}$ (note the process $X$ may be degenerate).}

\subsubsection*{Second set of hypotheses: Markovian assumptions} The
following assumptions will mostly be used in Section~\ref{section4MC}.
\HX, \HFp, and \HHp~give us a Markov representation for solutions
of the discrete BSDEs (see\vspace*{1pt} equation (\ref{eqmarkov}) later). \HGp~is
{used} for obtaining (model free) estimates on regression errors. We
also include additional optional assumptions, \HGt, on the terminal
condition to obtain tighter estimates on $Z_i$ (see Corollary~\ref
{coras} and subsequent remarks).
\begin{enumerate}[\HX]
\item[\HX] $X$ is a Markov chain in $\R^d$ $(1 \le d < +\infty)$
adapted to $(\mathcal{F}_{t_i})_i$.
For every $i<N$ and $j>i$, there exist $\mathcal{G}_{i}\otimes
\mathcal{B}(\R
^d)$-measurable functions $\Xmarkov ij\dvtx  \Omega\times\R^d
\rightarrow
\R^d$
where $\mathcal{G}_{i} \subset\mathcal{F}_T$ is independent of
$\mathcal{F}_{t_i}$, such that
$X_{j} = \Xmarkov{i}{j} (X_{i})$.
\item[\HGp]
\begin{enumerate}[(ii)]
\item[(i)] $\xi$ is a bounded $\mathcal{F}_T$-measurable random variable:
$C_\xi:=\P-\operatorname{ess}\sup_{\omega}\llvert\xi(\omega)\rrvert
<+\infty$.
\item[(ii)] $\xi$ is of form $\xi:=\Phi(X_N)$ for a {bounded},
measurable function $\Phi$.
\end{enumerate}
\item[\HGt] In addition to \HGp, for some $\thetaP\in[0,1]$
and a finite constant $C_\Phi\geq0$,
we have $\llvert\xi-\E_i\xi\rrvert_{2,i}\leq C_\Phi
(T-t_i)^{\thetaP/2}$ for any
$i\in\{0,\ldots,N\}$.
\item[\HFp] For every $i<N$, the driver is of the form $f_i(\omega
,y,z)=f_i(X_i(\omega),y,z)$, and
$(x,y,z)\mapsto f_i(x,y,z)$ is $\mathcal{B}(\R^d)\otimes\mathcal
{B}(\R) \otimes
\mathcal{B}
((\R^{q})^\top)$-measurable and \HF~is satisfied.
\item[\HHp] In addition to \HH, for every $i<N$ and $j>i$,
there is a function $\MWmarkov{i}{j}\dvtx  \Omega\times\R^d \rightarrow
(\R^{q})^\top$ that is $\mathcal{G}_{i}\otimes\mathcal{B}(\R
^d)$-measurable, where
$\mathcal{G}_{i} \subset\mathcal{F}_T$ is independent of $\mathcal
{F}_{t_i}$,
such that $H^{(i)}_{j}: = \MWmarkov i{j}(X_{i})$.
\end{enumerate}
\emph{Comments.} \HX~is usually satisfied when $X$ is the solution of
SDE or its Euler scheme built on the time-grid $\pi$.

The statement of \HGt~is inspired by the fractional smoothness
condition of \cite{gobemakh10}, although somewhat stronger.

It is satisfied,
for instance, if the terminal function $\Phi$ is $\thetaP$-H\"older
continuous and if the Markov chain satisfies $\mathbb{E}_i[\llvert
X_N -
X_i\rrvert^2] \le
C_X (T-t_i)$. This is {a reasonable} assumption on the Markov chain,
since it is satisfied by a diffusion process (possibly including
bounded jumps) with bounded coefficients and also by its Euler approximation.
Indeed, we have
%
\begin{equation}
\nonumber
\llvert\xi- \mathbb{E}_i\xi\rrvert
_{2,i}\leq\bigl\llvert\Phi(X_N) - \Phi(X_i)
\bigr\rrvert_{2,i}\leq C_\Phi\bigl(C_X
(T-t_i)\bigr)^{{ \thetaP}/{2 }}.
\end{equation}

\HHp~is satisfied by the Malliavin weights (\ref
{eq4malweights}) under various 
conditions. It is valid for the example $X_t=g(t,W_t)$ mentioned before.
Consider now the more complex case of the SDE with (deterministic)
coefficients $b(t,x)$ for the drift and $(\sigma_1(t,x),\ldots,\sigma
_d(t,x))$ for the diffusion ($q=d$) both having first space-derivatives
that are uniformly bounded.
Recall the relation for the Malliavin derivative of a SDE given by
\begin{eqnarray*}
D_tX_r= \nabla X_r \nabla
X_t^{-1} \sigma(t,X_t)\1_{t \le r}=
\nabla_x X^{t,x}_r\mid_{x=X_t}
\sigma(t,X_t)\1_{t \le r},
\end{eqnarray*}
where $X^{t,x}$ denotes the SDE solution starting from $x$ at time $t$,
and $\nabla X_t:=\nabla_x X_t^{0,x}$ for $\nabla_x X_t^{0,x}$ solving
the $(d \times d)$-dimensional, matrix valued linear SDE
\[
\nabla_x X^{t,x}_{r} = I_d + \int
_{t}^r 
\nabla_x b\bigl(u,\XM
tx_{u}\bigr) \nabla_x X^{ t,x}_u \,\mathrm{d}u +
\sum_{j=1}^d \int_{t}^r
\nabla_x\sigma_j\bigl( u,\XM
tx_{u}\bigr) \nabla_x X^{ t,x}_u
\,\mathrm{d}W_{j,u}.
\]
Then, it is an easy exercise to prove that if $\sigma$ and its inverse
are uniformly bounded, then \HHp~is fulfilled.

\section{Stability}\label{sectionstability}

\subsection{Gronwall type inequalities}\label{sectiongronwall}
Here we gather deterministic inequalities frequently used throughout
the paper.
{These inequalities are crucial due to novel technical problems caused
by the Malliavin weights. }
They show how linear inequalities with singular coefficients propagate.
They take the form of unusual Gronwall type inequalities.
Their proofs are postponed to Appendix~\ref{appendixprooflemintegrals}.
We assume that $\pi$ is in the class of time-grids satisfying \HFc.

\begin{lemma}
\label{lemintegrals}
Let
{$\alpha,\beta>0 $ be finite. }
There exists a finite constant $\intB{\alpha}{\beta}\geq0$ depending
only on $R_{\pi}$, $\alpha$ and $\beta$ (but not on the time-grid)
such that,
for any $0\le i < k \le N$,
\begin{eqnarray*}
\sum_{j=i}^{k-1}{\Delta_j \over(t_{k} - t_{j})^{1 - \alpha} }
&\le&\intB{\alpha} {1}(t_{k} - t_{i})^{\alpha},
\\
\sum_{j=i+1}^{k-1}{\Delta_{j} \over(t_{k} - t_{j})^{1 -
\alpha
}(t_{j} - t_{i})^{1 - \beta} }
&\le&\intB{\alpha} {\beta}(t_{k} - t_{i})^{\alpha+ \beta-1}.
\end{eqnarray*}
\end{lemma}

\begin{lemma}[(Exponent improvement in recursive equations)]
\label{lemiterationgen1}
Let $\alpha\geq0, \beta\in(0,\frac{ 1 }{2 }]$ and $k\in\{
0,\ldots,N-1\}$.
Suppose that, for a {finite} constant $C_{u}\geq0$, the finite
{non-negative} real-valued sequences $\{u_l\}_{l \ge k}$
and $\{w_l\}_{l\ge k}$ satisfy
%
\begin{eqnarray}
u_j \le w_j +C_{u}\sum
_{l=j+1}^{N-1} \frac{u_l \Delta
_l}{(T-t_l)^{{{1}/2 - \beta}}(t_l - t_j)^{{{1}/2-\alpha}}},
\qquad k\leq j \leq
N. \label{eqiterationfeed}
\end{eqnarray}
Then, for two {finite} constants {$\cw\geq0$ and $\cuh\geq0$} that
depend only on $C_{u}, T, \alpha, \beta$ and $R_\pi$,
%
\begin{eqnarray}\label{eqiterationfeed2}
u_j &\le& \cw w_j + \cw\sum
_{l=j+1}^{N-1} \frac{w_l \Delta
_l}{(T-t_l)^{{{1}/2 - \beta}}(t_l - t_j)^{{{1}/2-\alpha}}}
\nonumber\\[-8pt]\\[-8pt]\nonumber
&&{} +\cuh\sum
_{l=j+1}^{N-1} \frac{u_l \Delta_l}{(T-t_l)^{{{1}/2 -
\beta
}}}, \qquad k\leq j \leq
N.
\end{eqnarray}
%
\end{lemma}

\begin{lemma}[(Intermediate a priori estimates)]\label{lemiterationgen2}
Let $\alpha\geq0, \beta\in(0,\frac{ 1 }{2 }]$ and $k\in\{0,\ldots,N-1\}
$. Assume that the finite {non-negative} real-valued sequences $\{
u_l\}_{l \ge k}$
and $\{w_l\}_{l\ge k}$ satisfy (\ref{eqiterationfeed2}) for finite
constants $\cw\geq0$ and $\cuh\geq0$.
Then, for any finite $\gamma> 0$, there is a {finite} constant $\cwg
{\gamma}{\geq0}$ (depending only on $\cw$, $\cuh$, $T$, $\alpha$,
$\beta$, {$R_\pi$} and $\gamma$)
such that
%
\begin{eqnarray}\label{eqintuw}
&& \sum_{l=j+1}^{N-1} \frac{u_l \Delta_{l}}{(T-t_{l})^{1/2 - \beta
}(t_l -t_j)^{1-\gamma} }
\nonumber\\[-8pt]\\[-8pt]\nonumber
&&\quad  \le
\cwg{\gamma} \sum_{l=j+1}^{N-1}
\frac{w_l\Delta
_l}{(T-t_l)^{1/2
- \beta}(t_l -t_j)^{1-\gamma}}, \qquad k\leq j \leq N.
\end{eqnarray}
\end{lemma}

Plugging (\ref{eqintuw}) with $\gamma=\frac{1}2+\alpha$ into
(\ref{eqiterationfeed}) gives a ready-to-use result.

\begin{proposition}[(Final a priori estimates)]\label{coroiterationgen2}
Under the assumptions of Lemma~\ref{lemiterationgen1}, (\ref
{eqiterationfeed}) implies
\begin{eqnarray*}
u_j \le w_j +\cwg{{1}/2+
\alpha}C_{u}\sum_{l=j+1}^{N-1}
\frac{w_l
\Delta_l}{(T-t_l)^{{{1}/2 - \beta}}(t_l - t_j)^{{{1}/2-\alpha
}}}, \qquad k\leq j \leq N.
\end{eqnarray*}
\end{proposition}
%

\subsection{Stability of discrete BSDEs with Malliavin weights}\label{sectionstab}

Suppose that $(Y_1,Z_1)$ (resp., $( Y_2, Z_2)$) solves a {MWDP}
with terminal condition/driver $(\xi_1, f_{1,i})$ (resp., $( \xi_2, f_{2,i})$).
We are interested in {obtaining estimates on} the differences $(Y_1 -
Y_2, Z_1 - Z_2)$.
{To give a notion of how stability estimates are used, the processes
$(Y_{1},Z_{1})$ are typically obtained by construction.
For example, in Section~\ref{sectionas}, they are $(0,0)$, whereas in
the proof of Theorem~\ref{thmMCerr}, they are a set of processes
determined from a series of arguments based on conditioning {w.r.t.}
the Monte Carlo samples.
One then applies the stability estimates based on a priori knowledge
that what stands on the right-hand side is beneficial to the computations.
In Corollary~\ref{coras}, for example, the right-hand side yields
almost sure bounds for the processes $(Y,Z)$.
}
{We note that the assumptions on the drivers in this section are
somewhat weaker than the general assumptions of Section~\ref
{section4assumptions}.}
The driver $f_{1,i}(y,z)$ {does not have} to be Lipschitz continuous, but
we assume that each $f_{1,i}( Y_{1,i+1},Z_{1,i})$ is in $\L_2(\mathcal{F}_T)$
so that $Y_{1,i}$ and $Z_{1,i}$ are also square integrable for any $i$
(thanks to \HH). The driver $ f_{2,i}(y,z)$ is locally Lipschitz
continuous w.r.t. $(y,z)$ as in \HFa, which
is crucial for the validity of the a priori estimates.
Additionally, we do not insist that the drivers be adapted, which will
be needed in the setting of sample-dependant drivers.
Define
\begin{eqnarray*}
\Delta Y&:=& Y_1  - Y_2, \qquad\Delta Z:=
Z_1 - Z_2, \qquad\Delta\xi:= \xi_1 -
\xi_2,
\\
\Delta f_i&:=& f_{1,i}(Y_{1,{i+1}},Z_{1,i})
- f_{2,i}({Y_{1,{i+1}},Z_{1,i}}).
\end{eqnarray*}
{Let $k\in\{0,\ldots,N-1\}$ be fixed: throughout this subsection,
$\mathcal{F}
_{t_k}$-conditional $\L_2$-norms are considered and we recall the
notation $\llvert U\rrvert_{2,k}:= \sqrt{\E_k[\llvert
U\rrvert^2]}$} for any square integrable
random variable~$U$. For $j\ge{k}$, define
\[
\llvert{\Theta_j}\rrvert_{2,k}:= \llvert\Delta
Y_{{j+1}}\rrvert_{2,k} + \llvert\Delta Z_j
\rrvert_{2,k}.
\]
Using \HH, we obtain $\mathbb{E}_i[\Delta\xi\MW iN] = \mathbb
{E}_i[(\Delta\xi- \mathbb{E}_i
\Delta\xi)\MW iN]$ {and
%
\begin{eqnarray}\label{eqboundmalliavin}
\bigl\llvert\mathbb{E}_i\bigl[\Delta\xi\MW iN\bigr]\bigr
\rrvert^2 &\leq& \mathbb{E}_i\bigl[\llvert\Delta\xi-
\mathbb{E}_i \Delta\xi\rrvert^2\bigr]
\frac{ C_{M}^2
}{(t_N-t_i) },
\nonumber\\[-8pt]\\[-8pt]\nonumber
\bigl\llvert\E_i z\bigl[ \Delta
f_j \MW ij \bigr]\bigr\rrvert^2 &\le&\frac
{C_M^2 \E_i[\llvert\Delta
f_j\rrvert^2]}{t_j - t_i},
\qquad j \ge i+1.
\end{eqnarray}
Combining this kind of estimates} with \HFa~and the triangle
inequality, our stability equations {(for $k\leq i$)} are
%
\begin{eqnarray}
\llvert\Delta Y_i\rrvert_{2,k} & \le&\llvert\Delta\xi
\rrvert_{2,k} + \sum_{j = i}^{N-1}
\llvert\Delta f_j\rrvert_{2,k} \Delta_j +
\sum_{j=i}^{N-1} \frac{L_{f_2}\llvert\Theta
_j\rrvert_{2,k}}{(T-t_j)^{(1-\thetaL
)/2}}
\Delta_j, \label{eqbounddelY2}
\\
\llvert\Delta Z_i\rrvert_{2,k} & \le&\frac{C_{M}\llvert
\Delta\xi- {\mathbb{E}_i\Delta
\xi
}\rrvert_{2,k}}{\sqrt{T-t_i}}+ \sum_{j = i+1}^{N-1} \frac{C_{M}\llvert\Delta f_j\rrvert
_{2,k}}{\sqrt{t_j - t_i}}
\Delta_j
\nonumber\\[-8pt]\label{eqbounddelZ2} \\[-8pt]\nonumber
&&{}  + \sum_{j=i+1}^{N-1}
\frac{L_{f_2}C_{M}\llvert\Theta
_j\rrvert_{2,k}}{(T-t_j)^{(1-\thetaL)/2}\sqrt{t_j - t_i}} \Delta_j.
\end{eqnarray}

\begin{proposition}
\label{propstability}
Taking $\alpha=0 $, $\beta= \thetaL/2 $ and
$C_{u}=L_{f_2}(C_{M}+\sqrt
T)$ in Lemmas~\ref{lemiterationgen1} and~\ref{lemiterationgen2},
recall the constant $\cwg{\gamma}$. Assume that $\xi_j$ is in $\L
_2(\mathcal{F}_T)$.
Moreover, for each $i \in\{0, \ldots, N-1\}$, assume that
$f_{1,i}(Y_{1,i+1},Z_{1,i})$ is in $\L_2(\mathcal{F}_T)$ and
$ f_{2,i}(y,z)$ is
locally Lipschitz continuous w.r.t. $y$ and $z$ as in \HFa, with a
constant $L_{ f_{2}}$.
Then, under \HH, we have
%
\begin{eqnarray*}
\llvert\Delta Y_i\rrvert_{2,k}& \le&\aprioriConst{1}
{y}\llvert\Delta\xi\rrvert_{2,k} +\aprioriConst{2} {y}\sum
_{{j=i}}^{N-1} \llvert\Delta f_j\rrvert
_{2,k} \Delta_j,\qquad0\leq k \leq i \leq N,
\\
\llvert\Delta Z_i\rrvert_{2,k} & \le&\aprioriConst{1}
{z}\frac{ \llvert\Delta\xi- \mathbb
{E}_i\Delta\xi
\rrvert_{2,k}}{\sqrt{T-t_i}} + \aprioriConst{2} {z} \sum_{j = i+1}^{N-1}
\frac{\llvert\Delta
f_j\rrvert_{2,k}}{\sqrt{t_j - t_i}} \Delta_j
\\
&&{}+ \aprioriConst{3} {z}
\llvert\Delta
\xi\rrvert_{2,k} (T-t_{i})^{{ \thetaL}/{2
}},\qquad0\leq k
\leq i < N,
\end{eqnarray*}
where the above constants 
can be written explicitly:
%
\begin{eqnarray*}
\aprioriConst{1} {y} &:=& 1+ L_{f_2} \cwg{1}\biggl( C_{M}
\intB{\thetaL/2} {1} + \intB{{1}/2 +\thetaL/2} {1} \sqrt T
\biggr) T^{\thetaL/2},
\\
\aprioriConst2y &:=& 1 + L_{f_2} \cwg{1} (C_{M}+ \sqrt T )
\intB{{\thetaL/2}} {1} {T^{\thetaL/2}},
\\
\aprioriConst1z &:=& C_{M}\biggl(1 + L_{f_2}\cwg{
{1}/2} {C_{M}} \intB{
{{\thetaL}}/2} {
{1}/2} {T^{\thetaL/2}}\biggr),
\\
\aprioriConst2z &:=& C_{M}\biggl(1 + L_{f_2}\cwg{
{1}/2} (C_{M}+ \sqrt T) \intB{\thetaL/2} {
{1}/2} {T^{\thetaL/2}}\biggr),
\\
\aprioriConst3z &:=& {C_{M}} L_{f_2}\cwg{
{1}/2} \intB{{1}/2 + \thetaL/2} {
{1}/2}.
\end{eqnarray*}
\end{proposition}

\begin{pf}
Using (\ref{eqbounddelY2}) and (\ref{eqbounddelZ2}), we obtain
%
\begin{eqnarray}\label{eqboundtheta2}
\llvert\Theta_j\rrvert_{2,k} & \le&
C_{M}\frac{
\llvert\Delta\xi{- \mathbb
{E}_j\Delta\xi
} \rrvert_{2,k}}{\sqrt{T-t_j}} + {\llvert\Delta\xi\rrvert_{2,k} }
+(C_{M}+\sqrt T) \sum_{l = j+1}^{N-1}
\frac{\llvert\Delta
f_l\rrvert_{2,k}\Delta
_l}{\sqrt{t_l - t_j}}
\nonumber\\[-8pt]\\[-8pt]\nonumber
&&{}+ (C_{M}+\sqrt T) \sum_{l=j+1}^{N-1}
\frac{L_{f_2} \llvert
\Theta
_l\rrvert_{2,k}\Delta_l}{(T-t_l)^{(1-\thetaL)/2}\sqrt{t_l -
t_j}},\qquad j\geq k.
\end{eqnarray}
%
%

\textit{Upper bound for} (\ref{eqboundtheta2}).
We apply
Lemmas~\ref{lemiterationgen1} and~\ref{lemiterationgen2}
under the setting $u_j=\llvert\Theta_j\rrvert_{2,k}$, $w_j=
C_{M}\frac{ \llvert\Delta\xi- \mathbb{E}_j\Delta\xi\rrvert
_{2,k}}{\sqrt
{T-t_j}} +\llvert\Delta
\xi\rrvert_{2,k}+
(C_{M}+\sqrt T) \sum_{l = j+1}^{N-1} \frac{\llvert\Delta
f_l\rrvert_{2,k}\Delta
_l}{\sqrt{t_l - t_j}} $, $\alpha=0$, $\beta=\frac\thetaL2$,
$C_{u}=L_{f_2}(C_{M}+\sqrt T)$.
To make results fully explicit, we first need to upper bound quantities
of the form ($\gamma>0$)
\[
\mathcal{I}^{(\gamma)}_{j+1}:= \sum_{l=j+1}^{N-1}
\frac{w_l\Delta_l}{(T-t_l)^{1/2 -
\thetaL/
2}(t_l -t_j)^{1-\gamma}}.
\]
Using that
$\llvert\Delta\xi- \mathbb{E}_l\Delta\xi\rrvert_{2,k}$
is non-increasing in
$l$ and
Lemma~\ref{lemintegrals}, we obtain
%
\begin{eqnarray} \label{eqigamma}
\mathcal{I}^{(\gamma)}_{j+1}&=&\sum
_{l=j+1}^{N-1}
\Biggl( C_{M}\frac{
\llvert\Delta\xi- \mathbb{E}_l\Delta\xi\rrvert
_{2,k}}{\sqrt{T-t_l}} +\llvert\Delta
\xi\rrvert_{2,k}+
(C_{M}+\sqrt T) \sum_{r = l+1}^{N-1} \frac{\llvert\Delta
f_r\rrvert_{2,k}\Delta
_r}{\sqrt{t_r - t_l}} \Biggr)\Delta_l\nonumber
\\
&&\hspace*{24pt}{} \Big/
\bigl(
(T-t_l)^{1/2 - \sfrac\thetaL2}(t_l -t_j)^{1-\gamma}\bigr)
\nonumber\\[-8pt]\\[-8pt]\nonumber
&\leq& C_{M}B_{\sfrac\thetaL2,\gamma}\frac{ \llvert
\Delta\xi-
\mathbb{E}
_{j+1}\Delta\xi\rrvert_{2,k}}{(T-t_j)^{1-\sfrac\thetaL2-\gamma}}
+B_{\sfrac{1}2+\sfrac\thetaL2,\gamma}\frac{\llvert\Delta\xi
\rrvert_{2,k}}{(T-t_j)^{1/2-\sfrac\thetaL2-\gamma}}\nonumber
\\
&&{} +(C_{M}+\sqrt T)B_{ \sfrac\thetaL2,\gamma}\sum
_{l =
j+2}^{N-1}\frac{\llvert\Delta f_l\rrvert_{2,k}\Delta
_l}{(t_l-t_j)^{1-\sfrac
\thetaL
2-\gamma}}.\nonumber
\end{eqnarray}
%

\textit{Upper bound for} $|\Delta Y_i|_{2,k}$.
Starting
from (\ref{eqbounddelY2}) and applying Lemma~\ref
{lemiterationgen2}, we get
\begin{eqnarray*}
\llvert\Delta Y_i\rrvert_{2,k} & \le&\llvert\Delta\xi
\rrvert_{2,k} + \sum_{j = i}^{N-1}
\llvert\Delta f_j\rrvert_{2,k} \Delta_j
+L_{f_2} \cwg{1} \mathcal{I}^{(1)}_{i};
\end{eqnarray*}
then using the estimate (\ref{eqigamma}) and $ \llvert\Delta
\xi-\E_i
\Delta\xi\rrvert_{2,k}\leq\llvert\Delta\xi\rrvert
_{2,k}$, we obtain the announced
inequality.

\textit{Upper bound of} $\llvert\Delta Z_i\rrvert_{2,k}$.
Starting from
(\ref{eqbounddelZ2}) and applying Lemma~\ref{lemiterationgen2},
we have
\begin{eqnarray*}
\llvert\Delta Z_i\rrvert_{2,k} & \le&\frac{C_{M}\llvert
\Delta\xi- \mathbb{E}_i\Delta
\xi
\rrvert_{2,k}}{\sqrt{T-t_i}}
+ \sum_{j = i+1}^{N-1} \frac{C_{M}\llvert\Delta f_j\rrvert
_{2,k}}{\sqrt{t_j - t_i}}
\Delta_j
\\
&&{}+ L_{f_2}C_{M}\cwg{\sfrac{1}2}
\mathcal{I}^{(\sfrac{ 1 }{2 })}_{i+1}; 
\end{eqnarray*}
therefore using the estimate (\ref{eqigamma}), we derive the
advertised upper bound on $\llvert\Delta Z_i\rrvert_{2,k}$.
\end{pf}

\subsection{Almost sure bounds}\label{sectionas}
{In order to obtain error estimates for the Monte Carlo scheme, we use
the model-free estimates of Proposition~\ref{propeqyzerrdecoM}.
Typically, these estimates require that the object one is trying to
approximate is bounded.
Therefore, the following almost sure bounds are crucial.}

\begin{corollary}
\label{coras}
Assume \HGpa, \HF~and \HH~and recall the constants $\aprioriConst\cdot
y$ and $\aprioriConst
\cdot z$ from Proposition~\ref{propstability} where $L_{f_2}$ is
replaced by $L_f$. Then, we have
%
\begin{eqnarray}
\llvert Y_i\rrvert& \le& C_{y,i}:= \aprioriConst1{y}
C_\xi+ \aprioriConst{2} {y}C_f B_{\thetaC,1}
(T-t_i)^{\thetaC}, \label{eqasy}
\\
\llvert Z_i\rrvert& \le& C_{z,i}:= \aprioriConst1z
\frac
{\operatorname{ess}\sup_{\omega}
\llvert\xi
- \mathbb{E}_i\xi\rrvert_{2,i}}{\sqrt{T-t_i}} + \frac{\aprioriConst
2{z}C_f B_{\thetaC,\sfrac{1}2}}{(T-t_i)^{1/2 -
\thetaC}} +\aprioriConst3{z} C_\xi(T-t_{i})^{{ \thetaL}/{2 }}.\qquad
\label{eqasz}
\end{eqnarray}
%
\end{corollary}

The above upper bounds a valid for terminal values $\xi$ admitted by
\HGpa, which is quite general.
Without any further information on $\xi$, we can derive the simple bounds
%
\begin{equation}
\label{equniformbound} \llvert Y_i\rrvert+\sqrt{T-t_i}
\llvert Z_i\rrvert\leq\cyz
\end{equation}
for an explicit, time uniform constant $\cyz$.
It may, however, be useful to take advantage of additional information
on $\xi$, {to obtain finer estimates on $C_{y,i}$ and $C_{z,i}$ with
the aim
of better tuning the parameters of the MWLS method (see Section~\ref
{sectioncomplexity}).}
Two situations are of particular interest.
\begin{itemize}
\item For zero terminal condition, $Y$ and $Z$ get smaller and smaller
as $t_i$ goes to $T$ as expected:
$\llvert Y_i\rrvert+\sqrt{T-t_i}\llvert Z_i\rrvert
\leq C(T-t_i)^{\thetaC}$ for a constant $C$
depending only on $\aprioriConst2y$, $\aprioriConst2z$, $C_f$,
$\thetaC$ and $R_\pi$.
This result is useful for variance reduction methods like the proxy
method of \cite{gobeturk14}, Section~2.2, the method of Martingale basis
\cite{bendstei12}, and the multilevel method of \cite{turk13}.
\item Under \HGt, we have $\llvert\xi-\E_i\xi\rrvert
_{2,i}\leq C_\Phi
(T-t_i)^{\thetaP/2}$, which leads to an improved estimate for $Z$:
$\llvert Z_i\rrvert\le C (T-t_i)^{-\sfrac{1}2 + \thetaC\land
(\sfrac{\thetaP}2)} $ for
some constant $C$ depending only on $\aprioriConst1z$, $\aprioriConst2z$,
$\aprioriConst3z$, $C_f$, $\thetaC$, $R_\pi$, $T$, $C_\xi$ and
$C_\Phi$.
\end{itemize}
This is why in the subsequent analysis, we keep track on the general
dependence on $i$ of the constants $C_{y,i}$ and $C_{z,i}$.

\begin{pf*}{Proof of Corollary~\ref{coras}}
$(0,0)$ is the solution of the {MWDP} with data $(\xi_1 \equiv0,
f_{1,i} \equiv0)$.
Applying Proposition~\ref{propstability} with $(Y_1,Z_1) = (0,0)$ and
$(Y_2,Z_2)=(Y,Z)$ yields
\begin{eqnarray}
\llvert Y_i\rrvert_{2,k} & \le& C^{(1)}_{y}
\llvert\xi\rrvert_{2,k} + C^{(2)}_{y}\sum
_{j=i}^{N-1}\bigl\llvert f_j(0,0)
\bigr\rrvert_{2,k} \Delta_j,
\nonumber
\\
\llvert Z_i\rrvert_{2,k} & \le&\frac{ C^{(1)}_{z}\llvert
\xi- \mathbb{E}_i\xi
\rrvert_{2,k}}{\sqrt{T-t_i}} +
C^{(2)}_{z} \sum_{j=i+1}^{N-1}
\frac{\llvert f_j(0,0)\rrvert_{2,k}}{\sqrt
{t_j-t_i}} \Delta_j+ C^{(3)}_{z}
\llvert\xi\rrvert_{2,k}(T-t_{i})^{\sfrac{ \thetaL
}{2 }},
\nonumber
\end{eqnarray}
for $i=0,\ldots,N-1$.
Taking $k=i$, plugging in the almost sure bounds on $\llvert\xi
\rrvert$
from \HGpa~and $\llvert f_j(0,0)\rrvert$ from \HFb, and
using Lemma~\ref{lemintegrals}
then yields the result.
\end{pf*}

\section{Monte Carlo regression scheme}\label{section4MC}

Throughout this section, the Markovian assumptions \HX, \HGp, \HFp~and
\HHp~are in force. The notation and preliminary results used in
this section overlap with \cite{gobeturk14}, Section~4, and we recall
and adapt them to the setting of MWLS in Section~\ref
{sectionMCprelims} for completeness.

\subsection{Preliminaries}\label{sectionMCprelims}
{An elegant property of} the Markovian assumptions is there are
measurable, deterministic (but unknown) functions $y_i(\cdot)\dvtx  \R^d
\rightarrow\R$ and
$z_i(\cdot)\dvtx  \R^d \rightarrow(\R^q)^\top$ for each $i \in\{
0,\ldots,N-1\}$
such that the solution $(Y_i,Z_i)_{0\le i \le N-1}$ of the discrete
BSDE (\ref{eqmalscheme}) is given by
%
\begin{equation}
(Y_{i},Z_i):= \bigl(y_i(X_{i}),z_i(X_{i})
\bigr). \label{eqmarkov}
\end{equation}
{In this section, we estimate these functions.} 
One needs to apply Lemma~\ref{lemcdnexp}
{below} combined with
$\mathcal{G}= \mathcal{G}_i$ -- defined in the assumptions \HX~and
\HHp~-- $U =
X_i$, and
\begin{eqnarray*}
F(x) &:=& \Phi\bigl(\Xmarkov iN(x)\bigr) + \sum_{k=i}^{N-1}
f_k \bigl(\Xmarkov ik(x),y_{k+1}\bigl(\Xmarkov i{k+1}
(x)\bigr), z_k\bigl(\Xmarkov ik (x)\bigr) \bigr)\Delta
_k\qquad\mbox{for } y_i(\cdot),
\end{eqnarray*}
and
\begin{eqnarray*}
F(x) &:=& \Phi\bigl(\Xmarkov iN(x)\bigr) \MWmarkov iN(x)
\\
&&{} + \sum
_{k=i+1}^{N-1} f_k \bigl(\Xmarkov
ik(x), y_{k+1}\bigl(\Xmarkov i{k+1} (x)\bigr), z_k\bigl(
\Xmarkov ik (x)\bigr) \bigr) \MWmarkov ik(x)\Delta_k\qquad\mbox{for
} z_i(\cdot).
\end{eqnarray*}

\begin{lemma}[(\cite{gobeturk14}, Lemma 4.1)]
\label{lemcdnexp}
Suppose that $\mathcal{G}$ and $\mathcal{H}$ are independent
sub-$\sigma
$-algebras of
$\mathcal{F}$.
{For $l \ge1$,} let $F\dvtx  \Omega\times\R^d \rightarrow{\R^l}$ be
bounded and $\mathcal{G}\otimes\mathcal{B}(\R^d) $-measurable,
and $U\dvtx  \Omega\rightarrow\R^d$ be $\mathcal{H}$-measurable.
Then,
$ \E[F(U) \mid\mathcal{H}] = j({U})$ where $j(h) = \E[F(h)]$ for
all $h \in
\R^d$.
\end{lemma}

Least-squares regression has its traditional implementation in
nonparametric statistics and signal processing \cite
{gyorkohlkrzywalk02}. In the traditional setting, the random object
is a pair of random variables {$(O,R)$ termed the ``observation'' $O$
and the ``response'' $R$. $R$ is considered to be some function of
$O$,} with the possible addition of noise, and one needs recover this function.
There are three important differences in the use of least-squares
regression methods in our setting, and for this reason we give a
definition of (ordinary) least-squares regression (OLS) that enables us
{to }approach our problems.
First, the {response} we consider is a nonlinear transformation of the
\emph{paths} of the Markov chain $X$ and the Malliavin weights $H$. We
are able {to simulate observations and responses (active learning)}
and we know the nonlinear function; what is unknown is the regression
function
, that is, the conditional expectation.
Therefore, OLS is defined in a way that easily enables path-dependence
and joint laws by defining the path of the Markov chain and Malliavin
weights as a single random variable, $\cX$, with law $\nu$. Secondly,
since we are in a dynamical setting, least-squares regressions will be
computed using independent clouds of simulations on each point of the time-grid.
This causes a dependence on an additional source of randomness in the
observations, namely the cloud of simulations from the preceding computations.
Therefore, OLS is defined to depend on \emph{two} probability spaces:
one for the preceding clouds $(\tilde\Omega,\tilde\mathcal
{F},\tilde\P)$,
and one for the current cloud {distribution} $(\R^l,\mathcal{B}(\R
^l),\nu)$.
Finally, we will make use of both general probability measures
(associated to the joint-law of the Markov chain and Malliavin weights)
and empirical measures.
The use of simulations to generate the empirical measure creates
dependency issues that are avoided when using laws, whence we make two
distinct definitions depending on which measure is in use.
We recall the {general} notation of \cite{gobeturk14}, Section~4.1,
for ordinary least-squares regression problems.

\begin{definition}[(Ordinary least-squares regression)]\label{defls}
For $l,l'\geq1$ and for probability spaces $(\tilde\Omega,\tilde
\mathcal{F},\tilde\P)$ and $(\R^l,\mathcal{B}(\R^l),\nu)$,
let $S$ be a $\tilde\mathcal{F}\otimes\mathcal{B}(\R^l)$-measurable
$\R^{l'}$-valued function such that $S(\omega,\cdot) \in\L
_2(\mathcal{B}
(\R
^l),\nu)$ for $\tilde\P$-a.e. $\omega\in\tilde\Omega$,
and $\mathcal{K}_{}$ a linear vector subspace of $\L_2(\mathcal
{B}(\R
^l),\nu)$
spanned by deterministic $\R^{l'}$-valued functions $\{p_k(\cdot ), k\geq
1\}$.
The least-squares approximation of $S$ in the space $\mathcal{K}_{}$ with
respect to $\nu$ is the ($\tilde\P\times\nu$-a.e.) unique,
$\tilde\mathcal{F}\otimes\mathcal{B}(\R^l)$-measurable function
$S^\star$
given by
%
\begin{equation}
\label{eqmcmls} S^\star(\omega, \cdot):=\arg\inf_{\phi\in\mathcal{K}_{} }
\int\bigl\llvert\phi(x) - S(\omega,x) \bigr\rrvert^2 \nu(\mathrm{d}x).
\end{equation}
We say that $S^\star$ solves $\OLS(S,\mathcal{K}_{},\nu)$.

{On the other hand, suppose that $\nu_M=\frac{1}M \sum_{m=1}^M \delta
_{\cX^{(m)}}$ is a discrete probability measure
on $(\R^l,\mathcal{B}(\R^l))$,
where $\delta_x$ is the Dirac measure on $x$ and $\cX^{(1)},\ldots,\cX
^{(M)}\dvtx  \tilde\Omega\rightarrow\R^l$ are i.i.d. random variables.
For an $\tilde\mathcal{F}\otimes\mathcal{B}(\R^l)$-measurable $\R
^{l'}$-valued
function $S$ such that $\llvert S (\omega,\cX^{(m)}(\omega
)
)\rrvert<
\infty$ for any $m$ and $\tilde\P$-a.e. $\omega\in\tilde\Omega$,
the least-squares approximation of $S$ in the space $\mathcal{K}_{}$ with
respect to $\nu_M$ is the ($\tilde\P$-a.e.) unique,
$\tilde\mathcal{F}\otimes\mathcal{B}(\R^l)$-measurable function
$S^\star
$ given by
%
\begin{equation}
\label{eqmclsempi} S^\star(\omega, \cdot):= \arg\inf_{\phi\in\mathcal
{K}_{} }
\frac{1}M \sum_{m=1}^M \bigl
\llvert\phi\bigl(\cX^{(m)}(\omega) \bigr) - S \bigl(\omega,
\cX^{(m)}(\omega) \bigr) \bigr\rrvert^2.
\end{equation}
We say that $S^\star$ solves $\OLS(S,\mathcal{K}_{},\nu_M)$.
}
\end{definition}

From (\ref{eqmarkov}),
the {MWDP} (\ref{eqmalscheme}) can be reformulated in terms of
Definition~\ref{defls}: taking for $\mathcal{K}_{i}^{(l')}$ any dense
subset {in the $\R^{l'}$-valued functions belonging to $\L_2(\mathcal
{B}(\R^d), \P\circ(X_i)^{-1} )$,}
for each $i\in\{0,\ldots,N-1\}$,
%
\begin{eqnarray}\label{eqMDPlsdef}
\cases{ \displaystyle y_i(\cdot) \mbox{ solves } \OLS\bigl(
\Obsa{Y,i} {{\vecx{x}^{(i)}}}, \mathcal{K}_{i}^{(1)}, \nu_{i} \bigr),
\cr
\quad\displaystyle\mbox{for } \Obsa{Y,i} {
\vecx{x}^{(i)}}:= \Phi(x_N) + \sum
_{k=i}^{N-1} f_k \bigl(x_k,y_{k+1}(x_{k+1}),
z_k(x_k) \bigr)\Delta_k,
\cr
z_i(\cdot) \mbox{ solves } \OLS\bigl( \Obsa{Z,i} {
\vecx{h}^{(i)},\vecx{x}^{(i)}}, \mathcal{K}_{i}^{(q)}, \nu_{i} \bigr),
\cr
\quad\displaystyle\mbox{for } \Obsa{Z,i} {
\vecx{h}^{(i)},\vecx{x}^{(i)}}:= \Phi(x_N)h_N
+ \sum_{{k=i+1}}^{N-1} f_k
\bigl(x_k,y_{k+1}(x_{k+1}), z_k(x_k)
\bigr) h_k \Delta_k,}
\end{eqnarray}
\begin{eqnarray}
\nu_{i}&:=& \P\circ\bigl(\MW i{i+1},\ldots, {\MW
i{N}},X_i,\ldots,X_N\bigr)^{-1},
\nonumber\\[-8pt]\label{eqvecxhvecxh}\\[-8pt]\nonumber
\vecx{h}^{(i)}&:=&(h_{i+1},\ldots,{h_{N}})
\in\bigl(\bigl(\R^q\bigr)^\top\bigr)^{{N-i}},
\qquad\vecx{x}^{(i)}:= (x_i,\ldots,x_N) \in
\bigl(\R^d\bigr)^{N-i+1}.
\end{eqnarray}
%
However, the above least-squares regressions encounter two
computational problems:
\begin{longlist}[(\textbf{CP1})]
\item[(\textbf{CP1})] $\L_2(\mathcal{B}({\R^d}),\P\circ(X_i )^{-1})$
is usually
infinite dimensional;
\item[(\textbf{CP2})] the integrals of the $\OLS$~in (\ref{eqMDPlsdef})
are presumably computed using the untractable law of
\[
\bigl(\MW i{i+1},\ldots, {\MW i{N}},X_i,\ldots,X_N\bigr).
\]
\end{longlist}
Therefore, the functions $y_i(\cdot)$ and $z_i(\cdot)$ are to be
approximated on
finite-dimensional function spaces with the sample-based empirical
version of the law, as described in the next subsection.

\subsection{Algorithm}\label{sectionMCalg}
The first computational problem (\textbf{CP1}) is handled using a pre-selected
\emph{finite-dimensional vector spaces}.

\begin{definition}[(Finite-dimensional approximation spaces)]\label
{deffindimspace}
For
$i\in\{0,\ldots, N-1\}$, we are given two finite functional linear
spaces of dimension $K_{Y,i}$ and $K_{Z,i}$
\[
\cases{ \mathcal{K}_{Y,i}:= \operatorname{span}\bigl\{\pl{Y,i}1,\ldots, \pl{Y,i}
{K_{Y,i} } \bigr\}, &\quad for $\pl{Y,i}k\dvtx \R^d \rightarrow \R$ s.t. $\mathbb{E}\bigl[\bigl\llvert\plx{Y,i}k{X_i}\bigr\rrvert
^2\bigr] < +\infty$,
\vspace*{5pt}\cr
\mathcal{K}_{Z,i}:= \operatorname{span}\bigl\{
\pl{Z,i}1,\ldots, \pl{Z,i} {K_{Z,i} } \bigr\}, & \quad for $\pl{Z,i}k\dvtx
\R^d \rightarrow{\bigl(\R^q\bigr)^\top}$ s.t. $\mathbb{E}\bigl[\bigl\llvert\plx{Z,i}k{X_i}\bigr\rrvert
^2\bigr] < +\infty$.}
\]
The function $y_i(\cdot)$ (resp., $z_{i}(\cdot)$) will be approximated
in the linear space $\mathcal{K}_{Y,i}$ (resp., $\mathcal{K}_{Z,i}$).
The best approximation errors are defined by
%
\[
\mathcal{E}^Y_{\mathrm{App.},i}:= \sqrt{\inf_{ \phi\in\mathcal{K}_{Y,i}
}
\mathbb{E} \bigl[ \bigl\llvert\phi(X_i) - y_i(X_i)
\bigr\rrvert^2 \bigr]}, \qquad\mathcal{E}^Z_{\mathrm{App.},i}:= \sqrt
{\inf_{ \phi\in\mathcal{K}_{Z,i}
}\mathbb{E} \bigl[ \bigl\llvert
\phi(X_i) - z_i(X_i)\bigr\rrvert
^2 \bigr]}.
\]
\end{definition}

The second computational problem (\textbf{CP2}) is solved using the \emph{empirical} measure built from independent simulations with distribution
$\nu_i$.
{The number of simulations is large enough to avoid having
under-determined systems of equations to solve.}

\begin{definition}[(Simulations and empirical measures)]\label{defsimsandempi}
For $i\in\{0,\ldots,N-1\}$, generate $M_i\geq K_{Y,i} \vee
K_{Z,i}$ independent copies $\mathcal{C}_i:= \{( \MW{i,m}{ },\X
{i,m}{} )\dvt
m = 1,\ldots,M_i \}$
of $(\MW i{ },X^{(i)}):=(\MW i{i+1},\ldots, {\MW i{N}},X_i,\ldots,X_N)$:
$\mathcal{C}_i$ forms a \emph{cloud of simulations} used for the
regression at time~$i$.
Denote by $\nu_{i,M}$ the empirical probability measure of the
$\mathcal{C}
_i$-simulations,
that is,
%
\begin{equation}
\label{eqempimeasure} \nu_{i,M}:= \frac{1}{M_i} \sum
_{m=1}^{M_i} \delta_{(\MW
{i,m}{i+1},\ldots,{\MW{i,m}{N}}, \X{i,m}{i},\ldots,\X{i,m}N )}.
\end{equation}
Furthermore, we assume that the clouds of simulations $(\mathcal
{C}_i\dvt 0\leq i <
N)$ are independently generated. All these random variables are defined
on a probability space $(\Om M, \F M{}, \Prob M)$.
\end{definition}

Observe that allowing time-dependency in the number of simulations
$M_i$ and in the vector spaces $\mathcal{K}_{Y,i}$ and $\mathcal{K}_{Z,i}$
is coherent with our setting of time-dependent local Lipschitz driver.

Denoting by $(\Omega, \mathcal{F}, \P)$ the probability space supporting
$(\MW
{0}{ },\ldots,\MW{N-1}{ }, X)$, which serves as a generic element for
the clouds of simulations, the full probability space used to analyze
our algorithm is the product space $(\bar\Omega, \bar\mathcal{F},
\bar\P
)=(\Omega, \mathcal{F},\P) \otimes(\Om M, \F M{}, \Prob M)$.
By a slight abuse of notation, we write $\P$ (resp., $\E$) to mean
$\bar\P$ (resp., $\bar\E$) from now on.

In what follows, extensive {use} will be made of conditioning on the
clouds of simulations. This is much in the spirit of the proof of
\cite{gobeturk14}, Theorem 4.11, and the arguments are based on the
following definition of $\sigma$-algebras.

\begin{definition}
\label{defsimcnd}
Define the $\sigma$-algebras
\[
\F*i:= {\sigma(\mathcal{C}_{i+1},\ldots,\mathcal{C}_{N-1})},
\qquad\F{M}i:= \F*i \vee\sigma\bigl(\X{i,m} {i}\dvt  1 \le m \le M_i\bigr).
\]
For every $i\in\{0,\ldots,N-1\}$, let
$\E^M_i[\cdot] $ (resp., $\P^M_i$) with respect to $\F Mi$.
\end{definition}

We now come to the definition of the MWLS algorithm: this is merely
the finite-dimensional version of (\ref{eqMDPlsdef}) plus a soft
truncation of the solutions using the truncation function $\mathcal{T}_\cdot(\cdot )$
(defined in Section~\ref{sectionnotation}).

\begin{definition}[(MWLS algorithm)]
\label{defmwls} Set $\yM{N}{\cdot}:= \Phi(\cdot)$. For each
$i=N-1,N-2,\ldots,0$, set the random functions $\yM i\cdot$ and $\zM
i\cdot$ recursively as follows.
\begin{enumerate}
\item First, define $\zM i\cdot:= \mathcal{T}_{\Czi} ( \psi
^{(M)}_{Z,i}(\cdot)
)$ where $\Czi$ is the almost sure bound of Corollary~\ref{coras}
and where
%
\begin{equation}
\label{eqPsiMz}
\hspace*{-10pt}\cases{ \displaystyle\psi^{(M)}_{Z,i}(\cdot)
\mbox{ solves } \OLS\bigl( \ObsMa{Z,i} {\vecx{h}^{(i)},
\vecx{x}^{(i)}}, \mathcal{K}_{Z,i}, \nu_{i,M}\bigr)
\vspace*{3pt}\cr
\displaystyle\quad\mbox{for } \ObsMa{Z,i} {\vecx{h}^{(i)},\vecx{x}^{(i)}}:=
\Phi(x_N)h_N + \sum_{k=i+1}^{N-1}
f_k \bigl(x_k,\yM{k+1} {x_{k+1}}, \zM
k{x_k} \bigr) h_k \Delta_k,}
\end{equation}
where $\vecx{h}^{(i)}, \vecx{x}^{(i)}, \nu_{i,M}$ are defined in
(\ref{eqvecxhvecxh}) and (\ref{eqempimeasure}).

\item Second and similarly, define $\yM i\cdot:= \mathcal{T}_{\Cyi
}
(\psi^{(M)}
_{Y,i}(\cdot) )$ where
%
\begin{equation}
\label{eqPsiMy} \cases{ \displaystyle\psi^{(M)}_{Y,i}(\cdot)
\mbox{ solves }\OLS\bigl(\ObsMa{Y,i} {\vecx{x}^{(i)}},
\mathcal{K}_{Y,i}, \nu_{i,M}\bigr)
\vspace*{3pt}\cr
\displaystyle \quad\mbox{for } \ObsMa{Y,i}
{\vecx{x}^{(i)}}:= \Phi(x_N) + \sum
_{k=i}^{N-1} f_k \bigl(x_k,
\yM{k+1} {x_{k+1}}, \zM{k} {x_k} \bigr) \Delta_k.}
\end{equation}
\end{enumerate}
\end{definition}

Before performing the error analysis, we state the following uniform
(resp., conditional variance) bounds on the functions $\ObsM{Y,i}
{\cdot}$ (resp., the $l$th coordinate of $\ObsM{Z,i} {\MW{i,m}{}, \X
{i,m}{} }$ for each $m$ and $l$).
{These bounds are used repeatedly in Section~\ref{sectionerr} in
conjunction with Proposition~\ref{proplsregproperties} in order to
obtain estimates on the conditional variance of the regressions. }
The proof is postponed to Appendix~\ref{appendixprooflemmalemobsbounds}.

\begin{lemma}
\label{lemobsbounds}
For all $i \in\{0,\ldots,N-1\}$,
there are finite constants $\ObsBd{y,i}\geq0$ and $\ObsBd{z,i}\geq0$
such that
\begin{eqnarray*}
\bigl\llvert\ObsMa{Y,i} {\vecx{x}^{(i)}}\bigr\rrvert& \le&\ObsBd{y,i},
\qquad\forall\vecx{x}^{(i)},
\\
\sum_{l=1}^q\var\bigl[\ObsMa{l,Z,i} {\MW{i,m}
{}, \X{i,m} {}} \mid\F{M}i \bigr]& \le&\ObsBd{z,i}^2, \qquad\forall m \in
\{1,\ldots,M_i\}.
\end{eqnarray*}
We can write a precise time-dependency
of the constants $ \ObsBd{y,i}$ and $\ObsBd{z,i}$:
%
\begin{eqnarray}\label{eqboundObsBd}
\ObsBd{y,i} &:=&c_1 C_\xi+ c_2
C_f (T-t_i)^\thetaC,
\nonumber\\[-8pt]\\[-8pt]\nonumber
\ObsBd{z,i}&:=&c_3 {C_\xi(T-t_i)^{-1/2}
} + c_4 C_f (T-t_i)^{\thetaC-{ 1 }/{2 }},
\end{eqnarray}
where $(c_{j})_{1\leq j\leq4}$ depend only on $(L_f, C_{M}, q,
\aprioriConst1{y},\aprioriConst2{y},\aprioriConst1{z},\aprioriConst
2{z},\aprioriConst3{z},T,R_\pi,\thetaL,\thetaC)$ (computed explicitely
in the proof).
\end{lemma}

The above time-dependency is to be used to derive convergence rates for
the complexity analysis.

\subsection{Main result: Error analysis}\label{sectionerr}
We precise the random norms used to quantify the error of MWLS.

\begin{definition}
Let
$\varphi\dvtx  \Om M \times\R^d \rightarrow\R$ or $(\R^q)^\top$ be
$\F M{}
\otimes\mathcal{B}({\R^d})$-measurable.
For each $i\in\{0,\ldots,N-1\}$, define the random norms
\begin{eqnarray}
\llVert\varphi\rrVert_{i,\infty} ^2:= \int
_{\R^d} \bigl\llvert\varphi(x)\bigr\rrvert^2 \P
\circ X_i^{-1} (\mathrm{d}x), \qquad\llVert\varphi\rrVert
_{i,M} ^2:= \frac{1}{M_i} \sum
_{m=1}^{M_i} \bigl\llvert\varphi\bigl(\X{i,m}i\bigr)\bigr
\rrvert^2.
\nonumber
\end{eqnarray}
\end{definition}

The accuracy of the MWLS algorithm is measured as follows:
\begin{eqnarray*}
\bar\mathcal{E}(Y,M,i)&: =&\sqrt{\E\bigl[ \bigl\llVert\yM i\cdot-
y_i(\cdot) \bigr\rrVert_{i,\infty} ^2 \bigr]},
\qquad\bar\mathcal{E}(Z,M,i):= \sqrt{\E\bigl[\bigl\llVert\zM i\cdot-
z_i(\cdot) \bigr\rrVert_{i,\infty} ^2 \bigr]},
\\
\mathcal{E}(Y,M,i) &:=& \sqrt{\E\bigl[\bigl\llVert\yM i \cdot-
y_i(\cdot)\bigr\rrVert_{i,M} ^2 \bigr]},
\qquad\mathcal{E}(Z,M,i):= \sqrt{ \mathbb{E} \bigl[ \bigl\llVert\zM i
\cdot- z_i(\cdot)\bigr\rrVert_{i,M} ^2 \bigr]
}.
\end{eqnarray*}
In our analysis, we will have to switch from errors in true measure
$\bar\mathcal{E}(\cdots)$ to errors in empirical measure $\mathcal
{E}(\cdots)$, and
vice-versa: this is not trivial since $(\yM i\cdot, \zM i\cdot )$ and the
empirical norm $\llVert\cdot\rrVert_{i,M} $ depend on the same
sample. However,
the switch can be performed
using concentration-of-measure estimates uniformly on a class of
functions \cite{gyorkohlkrzywalk02}, Chapter~9. We directly state
the ready-to-use result, which is a straightforward adaptation of \cite{gobeturk14}, Proposition 4.10, to our context.

\begin{proposition}\label{propeqyzerrdecoM}
Recall Definition~\ref{deffindimspace} and the constants $\Cyi$
(resp., $\Czi$) from Corollary~\ref{coras}, and define the \emph
{interdependence errors}
%
\begin{eqnarray*}
\mathcal{E}^Y_{\mathrm{Dep.},i}&:=&\Cyi\sqrt{\frac{2028 (K_{Y,i}+1)
\log(3M_i)}{M_i}},
\\
\mathcal{E}^Z_{\mathrm{Dep.},i}&:=&\Czi\sqrt{\frac{2028
(K_{Z,i}+1) q \log(3M_i) }{M_i}}.
\end{eqnarray*}
For each $i \in\{0,\ldots, N-1\}$, we have
\begin{eqnarray*}
\bar\mathcal{E}(Y,M,i) \le\sqrt2 \mathcal{E}(Y,M,i) + \mathcal
{E}^Y_{\mathrm{Dep.},i},\qquad\bar\mathcal{E}(Z,M,i) \le\sqrt2
\mathcal{E}(Z,M,i) + \mathcal{E}^Z_{\mathrm{Dep.},i}.
\end{eqnarray*}
\end{proposition}

The aim is to determine a rate of convergence for $\mathcal
{E}(Y,M,k):=(\E
[\llVert
y_k - y^M_k\rrVert^2_{k,M}])^{ 1 /2 }$ and $\mathcal{E}
(Z,M,k):=(\E[\llVert z_k -
z^M_k\rrVert^2_{k,M}])^{ 1 /2 }$ using the local error
terms $(\mathcal{E}
(k))_k$ defined below.

\begin{theorem}[(Global error of the MWLS algorithm)]
\label{thmMCerr}
For $0\le k\le N-1$, define
%
\begin{eqnarray}\label{eqerr}
\mathcal{E}(k)&:=& \mathcal{E}^Y_{\mathrm{App.},k+1} + \ObsBd{y,k+1}
\sqrt{\frac{ K_{Y,k+1}}{M_{k+1} }} + \mathcal{E}^Z_{\mathrm{App.},k}
\nonumber\\[-8pt]\\[-8pt]\nonumber
&&{}  +
\ObsBd{z,k} \sqrt{\frac{K_{Z,k}} { M_k} }+L_f \bigl( \mathcal{E}
^Y_{\mathrm{Dep.},k+1} + \mathcal{E}^Z_{\mathrm{Dep.},k}
\bigr).
\end{eqnarray}
For every $k \in\{0,\ldots,N-1\}$,
%
\begin{eqnarray}
\bigl(\E\bigl[\bigl\llVert y_k - y^M_k
\bigr\rrVert_{k,M} ^2\bigr] \bigr)^{1/2} & \le&
\mathcal{E}^Y_{\mathrm{App.},k} +\ObsBd{y,k} \sqrt{\frac{ K_{Y,k}}{M_k}}
+ \aprioriConst My \sum_{j=k}^{N-1}
\frac{ \mathcal{E}(j) \Delta_j}{(T-t_j)^{(1-\thetaL)/2}},\label
{eqMCyerr}
\\
\bigl(\E\bigl[\bigl\llVert z_k - z^M_k
\bigr\rrVert_{k,M} ^2\bigr] \bigr)^{1/2} & \le&
\mathcal{E}^Z_{\mathrm{App.},k} +\ObsBd{z,k} \sqrt{\frac
{K_{Z,k}} {
M_k} }
\nonumber\\[-8pt]\label{eqMCzerr} \\[-8pt]\nonumber
&&{} + \aprioriConst Mz \sum_{j=k+1}^{N-1}
\frac{\mathcal{E}(j)\Delta
_j}{(T-t_j)^{(1-\thetaL)/2}\sqrt{t_j-t_k}},
\end{eqnarray}
where, recalling the constant $\cwg\gamma$ from Lemma~\ref
{lemiterationgen2} (with $\alpha= 0$, $\beta= \frac\thetaL2$,
$\gamma\in\{\frac{1}2,1\}$ and $C_u = L_f( \sqrt2 C_{M}+ 4 \sqrt T)$),
\begin{eqnarray*}
\aprioriConst My &:=& 2+ 4L_f\cwg{1}\bigl(1+\intB{\sfrac\thetaL2}
{1} T^{\sfrac
\thetaL2}(C_{M}+ 2\sqrt T)\bigr),
\\
\aprioriConst Mz &:=& C_{M}+ \sqrt2C_{M}L_f
\cwg{\sfrac{1}2}\bigl(1+\intB{\sfrac\thetaL2} {{1}/2}
T^{\thetaL/2}(C_{M}+ 2\sqrt T)\bigr).
\end{eqnarray*}
\end{theorem}

\emph{Discussion}.
Observe that owing to Proposition~\ref{propeqyzerrdecoM}, similar
estimates (with modified constants) are valid for $\bar\mathcal
{E}(Y,M,k)=(\E
[\llVert y_k - y^M_k\rrVert^2_{k,\infty}])^{ 1 /2 }$ and $\bar\mathcal{E}
(Z,M,k)=(\E[\llVert z_k - z^M_k\rrVert^2_{k,\infty
}])^{ 1 /2 }$. The global
error (\ref{eqMCyerr})--(\ref{eqMCzerr}) is a weighted time-average
of three different errors.
\begin{longlist}[(2)]
\item[(1)] The\vspace*{1pt} contributions $\mathcal{E}^\cdot_{\mathrm{App.},\cdot}$ are
the best
approximation errors using the vector spaces of functions: this
accuracy is achieved asymptotically with an infinite number of
simulations (take $M_k \to+\infty$ in our estimates).\vspace*{1pt}
\item[(2)] The contributions $\sqrt{\frac{K_{\cdot,\cdot}} { M_\cdot}}$ are
the usual statistical error terms: the larger the number of simulations
or the smaller the dimensions of the vector spaces, the better the
estimation error.
\item[(3)] The contributions $\mathcal{E}^{\cdot}_{\mathrm{Dep.},\cdot}$ are
related to
the interdependencies between regressions at different times. This is
intrinsic to the dynamic programming equation with $N$ nested empirical
regressions.
\end{longlist}
However, due to Proposition~\ref{propeqyzerrdecoM}, the latter
contributions are of same magnitude as the statistical error terms (up
to logarithmic factors). Therefore roughly speaking, the global error
is of order of the best approximation errors plus the statistical
errors, as if there were a single regression problem \cite{gyorkohlkrzywalk02}, Theorem~11.1. In this sense, these error bounds are
optimal: it is not possible to improve the above estimates with respect
to the convergence rates (but only possibly with respect to the
constants). An optimal tuning of parameters is proposed in Section~\ref
{sectioncomplexity}.

In comparison to \cite{gobeturk14}, where a different Monte Carlo
regression scheme is analyzed, the upper bound for the global error has
a similar shape, but with two important differences.
\begin{itemize}
\item\emph{Norm on $Z$.} {In \cite{gobeturk14}, one uses the time
averaged {squared }$\L_2$-norm $\sum_i \E[ {\llVert\cdot\rrVert_{i,M}
^2} ]
\Delta_i$ to estimate the error in $Z$, whereas here the norm used is
time-wise. This leads to more informative error bounds. This is an
advantage of the discrete BSDE with Malliavin weights against the MDP
of \cite{gobeturk14}.}
\item\emph{Time-dependency.} The MWDP yields better estimates on
$y(\cdot )$ and $z(\cdot )$ w.r.t. time in the local error estimates, which
allows better parameters tuning and therefore better convergence rates
(see Section~\ref{sectioncomplexity}).
\end{itemize}

\subsection{Proof of Theorem \texorpdfstring{\protect\ref{thmMCerr}}{3.10}}\label{sectionproofmainthm}
\subsubsection{Preliminary results}

The following proposition lists  standard tools from the theory of
regression. They will be used repeatedly in the proof of Theorem~\ref
{thmMCerr}.
This proposition was also used in \cite{gobeturk14}, and we refer the
reader to that paper for the proof.
The two first properties are of deterministic nature, the two last are
probabilistic.
Item (iv) is stated in high generality; this readily allows its
further use in other regression-based Monte Carlo algorithms.

\begin{proposition}[(\cite{gobeturk14}, Proposition~4.12)]
\label{proplsregproperties}
With the notation of Definition~\ref{defls},
suppose that $\mathcal{K}_{}$ is
finite-dimensional and spanned by the functions $\{p_1(\cdot ),\ldots,p_K(\cdot )\}$.
Let $S^\star$ solve $\OLS(S,\mathcal{K}_{},\nu)$ (resp., $\OLS
(S,\mathcal{K}_{},\nu_M)$), according to (\ref{eqmcmls}) (resp.,
(\ref
{eqmclsempi})).
The following properties are satisfied:
\begin{longlist}[(iii)]
\item[(i)] Linearity: the mapping $S\mapsto S^\star$ is linear.
\item[(ii)] {Stability} property: $\llVert S^\star\rrVert_{\L
_2(\mathcal{B}(\R^l),\mu)}\leq\llVert S\rrVert_{\L_2(\mathcal{B}(\R
^l),\mu)}$, where $\mu= \nu$ (resp., $\mu= \nu_M$).
\item[(iii)] Conditional expectation solution:
in the case of the discrete probability
measure $\nu_M$,
assume additionally that the {sub-}$\sigma$-algebra $\mathcal
{Q}\subset
\tilde
\mathcal{F}$ is such that $ (p_j(\cX^{(1)}),\ldots,p_j(\cX
^{(M)})
)$ is
$\mathcal{Q}$-measurable for every $j\in\{1,\ldots,K
\}$. Setting
$S_\mathcal{Q}(\cX^{(m)}):=\tilde\E[S(\cX^{(m)})\mid\mathcal
{Q}]$ for each $m
\in\{
1,\ldots,M\}$,
then $\tilde\E[S^\star\mid\mathcal{Q}]$ solves $\OLS(
S_\mathcal{Q},
\mathcal{K}_{},
\nu_M )$.
\item[(iv)] Bounded conditional variance: in the case of the
discrete probability measure $\nu_M$, suppose that $S(\omega,x)$ is
$\mathcal{G}\otimes\mathcal{B}(\R^l)$-measurable, for $\mathcal
{G}\subset\tilde
\mathcal{F}$
independent of $\sigma(\Samp1,\ldots, \Samp M)$,
there exists a Borel measurable function $g\dvtx  \R^l \rightarrow
\mathcal{E}$,
for some Euclidean space $\mathcal{E}$, such that the random variables
$\{
p_j(\Samp m)\dvt  m =1,\ldots,M, j = 1,\ldots, K\}$ are $\mathcal{H}:=
\sigma
(g(\Samp m)\dvt  m =1,\ldots, M)$-measurable, and there is a finite
constant $\sigma^2 \geq0$ that uniformly bounds the conditional
variances $\tilde\E[\llvert S(\Samp{m})-\tilde\E(S(\Samp
{m})\mid\mathcal{G}
\vee\mathcal{H}
)\llvert^2 \mid\mathcal{G}\vee\mathcal{H} ]\leq
\sigma^2$
$\tilde\P$-a.s. and for
all $m
\in\{1,\ldots,M\}$. Then
\[
\tilde\E\bigl[ \bigl\llVert S^\star(\cdot) - \tilde\E
\bigl[S^\star(\cdot) \mid \mathcal{G}\vee\mathcal{H}\bigr]
\bigr\rrVert_{\L_2(\mathcal{B}(\R^l),\nu_M)}^2 \mid \mathcal
{G} \vee
\mathcal{H} \bigr] \le\sigma^2 \frac{ K }{M }.
\]
\end{longlist}
\end{proposition}

\emph{Intermediate processes and local error terms.}
Another technique we borrow from \cite{gobeturk14} is to introduce
intermediate, \emph{fictional} regressions based on the true solutions:
one replaces the full $\L_2$ space for the approximation space and the
true measure for the empirical measure in (\ref{eqMDPlsdef}).

For each $k\in\{0,\ldots,N-1\}$, recall the functions $\Obs
{Y,k}{\vecx
{x}^{(i)}}$ and $\Obs{Z,k}{\vecx{h}^{(i)},\vecx{x}^{(i)}}$ from
(\ref{eqMDPlsdef}),
the linear spaces $\mathcal{K}_{Y,k}$ and $\mathcal{K}_{Z,k}$ from Definition
\ref{deffindimspace}, and the empirical measure
$\nu_{k,M}$ from (\ref{eqempimeasure}), and set
\begin{eqnarray*}
&\displaystyle\psi_{Y,k}(\cdot) \quad
\mbox{solves}\quad\OLS\bigl( \Obsa{Y,k} {\vecx{x}^{(i)}},
\mathcal{K}_{Y,k}, \nu_{k,M} \bigr),&
\\
&\displaystyle\psi_{Z,k}(\cdot) \quad\mbox{solves}\quad\OLS\bigl( \Obsa{Z,k} {\vecx
{h}^{(i)},\vecx{x}^{(i)}}, \mathcal{K}_{Z,k},
\nu_{k,M} \bigr).&
\end{eqnarray*}

Note that these regressions are not {numerically} accessible, because
they require knowledge of the true solution to be applied.
After a series of conditioning arguments, based on Lemma~\ref{lembary}
below, the fictional regressions will eventually allow the use of
the stability estimates of Section~\ref{sectionstab}, and (after a
somewhat complex application of the Gronwall inequalities of
Section~\ref{sectiongronwall}) this will yield final result.

From Lemma~\ref{lemcdnexp} and our Markovian assumptions, observe that
\[
\bigl(\E^M_k\bigl[\Obsa{Y,k} {\X{k,m} {} }\bigr],
\E^M_k\bigl[ \Obsa{Z,k} {\MW{k,m} {},\X{k,m} {} }\bigr] \bigr)
= \bigl( y_k\bigl(\X{k,m} {k}\bigr), z_k\bigl(\X{k,m} {k}\bigr) \bigr)
\]
for each $m\in\{1,\ldots,M_k \}$ where $ (y_k(\cdot),z_k (\cdot)
)$ are the unknown functions defined in (\ref{eqmarkov}). Proposition
\ref{proplsregproperties}(iii) implies the first statement of the
following lemma. The second statement results from a direct interchange
of $\inf$ and $\E$, and from the identical distribution of
$(X^{(k,m)}_k)$ for all~$m$.

\begin{lemma}
\label{lembary}
For each $k\in\{0,\ldots,N-1\}$,
\begin{eqnarray*}
& \displaystyle\E^M_k\bigl[
\psi_{Y,k}(\cdot)\bigr] \quad\mbox{solves}\quad\OLS\bigl(
y_k(\cdot ), \mathcal{K}_{Y,k}, \nu_{k,M} \bigr),&
\\
& \displaystyle \E^M_k\bigl[\psi_{Z,k}(\cdot)\bigr] \quad
\mbox{solves}\quad\OLS\bigl( z_k(\cdot ), \mathcal{K}_{Z,k},
\nu_{k,M} \bigr).&
\end{eqnarray*}
In addition, recalling the local error terms $\mathcal{E}^Y_{\mathrm
{App.},k}$ and $\mathcal{E}^Z_{\mathrm{App.},k}$
from Definition~\ref{deffindimspace},
\begin{eqnarray*}
 \mathbb{E} \bigl[ \bigl\llVert\E^M_k\bigl[
\psi_{Y,k}(\cdot) \bigr]-y_k(\cdot) \bigr\rrVert
_{k,M} ^2 \bigr]&=&\mathbb{E} \Bigl[ \inf
_{ \phi\in\mathcal{K}_{Y,k} } \bigl\llVert\phi(\cdot) - y_k(\cdot)\bigr
\rrVert_{k,M} ^2 \Bigr]\leq\bigl(\mathcal{E}
^Y_{\mathrm{App.},k}\bigr)^2,
\\
 \mathbb{E} \bigl[ \bigl\llVert\E^M_k\bigl[
\psi_{Z,k}(\cdot) \bigr] -z_k(\cdot)\bigr\rrVert
_{k,M} ^2 \bigr]&=&\mathbb{E} \Bigl[ \inf
_{ \phi\in\mathcal{K}_{Z,k} } \bigl\llVert\phi(\cdot) - z_k(\cdot)\bigr
\rrVert_{k,M} ^2 \Bigr]\leq\bigl(\mathcal{E}
^Z_{\mathrm{App.},k}\bigr)^2.
\end{eqnarray*}
\end{lemma}

\subsubsection{Proof of Theorem \texorpdfstring{\protect\ref{thmMCerr}}{3.10}}
\begin{longlist}
\item[\textit{Step} 1: \textit{decomposition of the error on $Y$.}]
Recall the soft truncation function $\mathcal{T}_L(x):= (-L\vee x_1
\wedge L,\ldots, -L \vee x_n \wedge L)$ {for $x\in\R^n$}.
From the almost sure bounds of Corollary~\ref{coras}, $\mathcal
{T}_{\Cyk}
(y_k) = y_k$. Then, the Lipschitz continuity of $\mathcal{T}_{\Cyk}$ yields
$ \llVert y_k(\cdot) - \yM{k}{\cdot} \rrVert_{k,M} $ is
less than or
equal to
$\llVert y_k(\cdot) - \psi^{(M)}_{Y,k}(\cdot) \rrVert
_{k,M} $.
Using the triangle inequality for the $\llVert\cdot\rrVert_{k,M}
$-norm, it
follows that
%
\begin{eqnarray}\label{eqMCy1}
&& \bigl\llVert y_k(\cdot)- \yM k\cdot\bigr\rrVert_{k,M}
\nonumber\\[-8pt]\\[-8pt]\nonumber
&&\quad  \le\bigl\llVert y_k(\cdot) - \E^M_k
\bigl[ \psi_{Y,k}(\cdot) \bigr] \bigr\rrVert_{k,M} + \bigl
\llVert\mathbb{E} ^{M}_k\bigl[ \psi_{Y,k}(
\cdot)\bigr] - \psi^{(M)}_{Y,k}(\cdot) \bigr\rrVert
_{k,M}.
\end{eqnarray}
Because $\ObsM{Y,k}{\cdot}$ depends on $\zM{k}{\cdot}$ computed with
the same cloud of simulations $\mathcal{C}_k$ as that used to define
the OLS
solution $\psi^{(M)}_{Y,k}(\cdot)$, it raises some interdependency issues
that we solve by {making} a small perturbation {to the intermediate
processes} as follows {(compare with (\ref{eqMDPlsdef}) and~(\ref{eqPsiMy}))}:
for $\vecx{x}^{(k)} = (x_k,\ldots, x_N)$, define
\begin{eqnarray*}
&\displaystyle \tObsMa{Y,k} {\vecx{x}^{(k)}} := \Phi(x_N)
+f_k \bigl(x_k, \yM{k+1} {x_{k+1}},
z_{k}(x_k) \bigr) \Delta_k+ \sum
_{i=k+1}^{N-1} f_i \bigl(x_i,
\yM{i+1} {x_{i+1}}, \zM{i} {x_i} \bigr) \Delta_i,&
\\
&\displaystyle\tilde\psi^{(M)}_{Y,k}(\cdot) \quad\mbox{solves}\quad\OLS
\bigl(\tObsMa{Y,k} {\vecx{x}^{(k)}}, \mathcal{K}_{Y,k},
\nu_{k,M}\bigr).&
\end{eqnarray*}
This perturbation is not needed for the $Z$-component, because $\ObsM
{Z,k}{\vecx{h}^{(k)},\vecx{x}^{(k)}}$ depends {only} on the subsequent
clouds of simulations $\{\mathcal{C}_{j},j\ge k+1\}$.
Applying the $\L_2$-norm $\llvert\cdot\rrvert_2$, the
triangle inequality in
(\ref{eqMCy1}), and the first part of Lemma~\ref{lembary} yields
%
\begin{eqnarray} \label{eqMCy2}
\mathcal{E}(Y,M,k) & \le&\mathcal{E}^Y_{\mathrm{App.},k} + \bigl
\llvert\bigl\llVert\mathbb{E}^{M}_k\bigl[ \tilde\psi
^{(M)}_{Y,k}(\cdot)-\psi_{Y,k}(\cdot)\bigr]\bigr
\rrVert_{k,M} \bigr\rrvert_2\nonumber
\\
&&{}  + \bigl\llvert\bigl\llVert
\tilde\psi^{(M)}_{Y,k}(\cdot)- \mathbb{E}^{M}_k
\bigl[ \tilde\psi^{(M)} _{Y,k}(\cdot)\bigr]\bigr\rrVert
_{k,M} \bigr\rrvert_2
\\
&&{}+ \bigl\llvert\bigl\llVert\tilde\psi^{(M)}_{Y,k}(
\cdot)-\psi^{(M)}_{Y,k}(\cdot) \bigr\rrVert_{k,M}
\bigr\rrvert_2.\nonumber
\end{eqnarray}
Let us handle each term in the above inequality separately.
\begin{longlist}[$\rhd$]
\item[$\rhd$] Term $\llvert\llVert\mathbb{E}^{M}_k[
\tilde\psi^{(M)} _{Y,k}(\cdot)-\psi_{Y,k}(\cdot)]\rrVert
_{k,M} \rrvert_2$.
Set
\[
\txiM{Y, k}x:=\E\bigl(\tObsMa{Y,k} {{X}^{(k)}}-\Obsa{Y,k}
{{X}^{(k)}}\mid X^{(k)}_k=x,\F M{}\bigr).
\]
Recalling that $\tObsM{Y,k}{\vecx{x}^{(k)}}-\Obs{Y,k}{\vecx{x}^{(k)}}$
is built only using the clouds $\{\mathcal{C}_{j},j\geq k+1\}$,
it follows from {Lemma~\ref{lemcdnexp}} that $\mathbb
{E}^{M}_k[\tObsM{Y,k}{\X
{k,m}{ } }-\Obs{Y,k}{\X{k,m}{ } }]$ is equal to
$\txiM{Y,k}{\X{k,m}{k}}$ for every $m \in\{1,\ldots, M_k\}$.
Then, using Proposition~\ref{proplsregproperties}(i)~and~(iii),
$\mathbb{E}^{M}_k[ \tilde\psi^{(M)}_{Y,k}(\cdot)-\psi_{Y,k}(\cdot
)]$ solves
$\OLS(\txiM
{Y,k}\cdot, \mathcal{K}_{Y,k}, \nu_{k,M})$.
By Proposition~\ref{proplsregproperties}(ii),
\[
\E\bigl[\bigl\llVert\mathbb{E}^{M}_k\bigl[ \tilde
\psi^{(M)}_{Y,k}(\cdot)-\psi_{Y,k}(\cdot)\bigr]\bigr
\rrVert_{k,M} ^2 \bigr] \leq\E\bigl[ \bigl\llVert\txiM{Y,k}
\cdot\bigr\rrVert_{k,M} ^2 \bigr] = \E\bigl[\bigl(\txiM{Y,k}
{X_k}\bigr)^2 \bigr], %
\]
where the final equality follows from the fact that $\txiM{Y,k}{\cdot}$
is generated only using the simulations in the clouds $\{\mathcal{C}_j\dvt  j
>k\}$ and $\{X_k,\X{k,1}{k},\ldots,\X{k,M_k}{k}\}$ are identically
distributed.
Defining
%
\begin{equation}
\label{eqximy} \xiM{Y,k}x:= \E\bigl[\ObsMa{Y,k} {X^{(k)}}-\Obsa{Y,k}
{X^{(k)}}\mid X^{(k)}_k=x,\F M{}\bigr],
\end{equation}
the triangle inequality yields
\begin{eqnarray*}
\bigl\llvert\txiM{Y,k} {X_k } \bigr\rrvert_2 & \le&\bigl\llvert
\tObsMa{Y,k} {{X}^{(k)}} - \ObsMa{Y,k} {{X}^{(k)}} \bigr\rrvert
_2 + \bigl\llvert\xiM{Y,k} {X_k} \bigr\rrvert_2
\\
& \le&\bigl\llvert f_k\bigl(X_k,\yM{k+1}
{X_{k+1}}, \zM k{X_{k} } \bigr) - f_k
\bigl(X_k,\yM{k+1} {X_{k+1}}, z_k(X_{k})
\bigr) \bigr\rrvert_2 \Delta_k
\\
&&{}+ \bigl\llvert\xiM{Y,k}
{X_k} \bigr\rrvert_2
\\
& \le&\frac{L_f\Delta_k}{(T-t_k)^{1/2-\sfrac\thetaL2}} \bar\mathcal
{E}(Z,M,k) + \bigl\llvert\xiM{Y,k}
{X_k}\bigr\rrvert_2.
\end{eqnarray*}
\item[$\rhd$] Term $\llvert\llVert\tilde\psi
^{(M)}_{Y,k}(\cdot)- \mathbb{E}^{M}_k[ \tilde\psi^{(M)}_{Y,k}(\cdot
)]\rrVert_{k,M} \rrvert_2$.
Since $ \tObsM{Y,k}{\cdot}{ }$ depends only on the clouds $\{\mathcal{C}
_{j},j\geq
k+1\}$ and is bounded from above by $\ObsBd{y,k} $ (like $\ObsM
{Y,k}{\cdot}{ }$, see Lemma~\ref{lemobsbounds}), it follows from
Proposition~\ref{proplsregproperties}(iv) that $\llvert
\llVert\tilde\psi^{(M)} _{Y,k}(\cdot)- \mathbb{E}^{M}_k[
\tilde\psi^{(M)}_{Y,k}(\cdot)]\rrVert_{k,M} \rrvert_2$ is
bounded from
above by $\ObsBd{y,k} \sqrt{K_{Y,k} /M_k}$. This is similar to the
statistical error term in usual regression theory.

\item[$\rhd$] Term $\llvert\llVert\tilde\psi
^{(M)}_{Y,k}(\cdot)-\psi^{(M)} _{Y,k}(\cdot) \rrVert_{k,M}
\rrvert_2$.
Owing to Proposition~\ref{proplsregproperties}(i)~and~(ii), $\llVert
\tilde\psi^{(M)}_{Y,k}(\cdot)-\psi^{(M)}_{Y,k}(\cdot) \rrVert
_{k,M} ^2$ is
bounded from
above by
$\llVert\tObsM{Y,k}\cdot-\ObsM{Y,k}\cdot\rrVert_{k,M}
^2$, which equals
\begin{eqnarray*}
&& \frac{\Delta_k^2}{M_k} \sum_{m=1}^{M_k} \bigl
\llvert f_k\bigl(\X{k,m}k,\yMa{k+1} {\X{k,m} {k+1}}, \zMa k{\X{k,m} {k}}
\bigr) - f_k\bigl(\X{k,m}k,\yMa{k+1} {\X{k,m} {k+1}}, z_k\bigl(
\X{k,m} {k}\bigr) \bigr) \bigr\rrvert^2
\\
&&\quad \le\frac{ L_f^2 \Delta_k^2 \llVert z_k(\cdot) - \zM k\cdot
\rrVert_{k,M} ^2
}{ (T-t_k)^{1-\thetaL} }.
\end{eqnarray*}
Collecting the bounds on the three terms, substituting them into (\ref
{eqMCy2})
{and applying Proposition~\ref{propeqyzerrdecoM} }
yields
%
\begin{eqnarray}\label{eqMCy3}
\mathcal{E}(Y,M,k) & \le&\mathcal{E}^Y_{\mathrm{App.},k} + \bigl\llvert
\xi^*_{Y,k}(X_k) \bigr\rrvert_2 +
\ObsBd{y,k} \sqrt{\frac{ K_{Y,k}}{M_k}}
\nonumber\\[-8pt]\\[-8pt]\nonumber
&&{}  + \frac{L_f \Delta_k}{
(T-t_k)^{1/2-\sfrac\thetaL2}} \bigl\{ (1 + \sqrt2) {
\mathcal{E}(Z,M,k)} + \mathcal{E}^Z_{\mathrm
{Dep.},k} \bigr\}.
\end{eqnarray}
\end{longlist}
\end{longlist}

\begin{longlist}
\item[\textit{Step} 2: \textit{decomposition of the error on $Z$.}]
Analogously to (\ref{eqMCy2}), one obtains the upper bound
\begin{eqnarray}
\mathcal{E}(Z,M,k) & \le&\mathcal{E}^Z_{\mathrm{App.},k} + \bigl
\llvert\bigl\llVert\mathbb{E}^{M}_k\bigl[
\psi^{(M)}_{Z,k}(\cdot)-\psi_{Z,k}(\cdot)\bigr]\bigr
\rrVert_{k,M} \bigr\rrvert_2 + \bigl\llvert\bigl\llVert
\psi^{(M)}_{Z,k}(\cdot)- \mathbb{E}^{M}_k
\bigl[ \psi^{(M)}_{Z,k}(\cdot)\bigr]\bigr\rrVert
_{k,M} \bigr\rrvert_2.
\nonumber
\end{eqnarray}
Since $\ObsM{Z,k}{\cdot} $ depends only on the clouds $\{\mathcal
{C}_{j},j\geq
k+1\}
$ and the $\F M k$-conditional variance of $\ObsM{Z,k}{\MW{k,m}{},\X
{k,m}{} }$ is bounded from above by $ \ObsBd{z,k}^2 $ for all $m$ (see
Lemma~\ref{lemobsbounds}), it follows from Proposition~\ref
{proplsregproperties}(iv) that $\llvert\llVert\psi
^{(M)}_{Z,k}(\cdot)- \mathbb{E} ^{M}_k[ \psi^{(M)}_{Z,k}(\cdot
)]\rrVert_{k,M} \rrvert_2$ is bounded from
above by
$\ObsBd
{z,k} \sqrt{K_{Z,k} / M_k}$.
Defining
%
\begin{equation}
\label{eqximz} \xi^*_{Z, k}(x):=\E\bigl[\ObsMa{Z,k} {\MW{k}
{},{X}^{(k)}}-\Obsa{Z,k} {\MW{k} {},{X}^{(k)}}\mid
X^{(k)}_k=x,\F M{}\bigr],
\end{equation}
it follows that
$\mathbb{E}^{M}_k[ \psi^{(M)}_{Z,k}(\cdot)-\psi_{Z,k}(\cdot)]$
solves $\OLS(\xi
^*_{Z,k}(\cdot), \mathcal{K}_{Z,k}, \nu_{k,M})$. Therefore,
%
\begin{eqnarray}
\mathcal{E}(Z,M,k) & \le&\mathcal{E}^Z_{\mathrm{App.},k} + \bigl\llvert
\xiM{Z,k} {X_k} \bigr\rrvert_2 +\ObsBd{z,k} \sqrt{
\frac{K_{Z,k}} { M_k} }. \label{eqMCz1}
\end{eqnarray}
\end{longlist}

\begin{longlist}
\item[\textit{Step} 3: \textit{error propagation and a priori estimates.}]
Observe that $(\xiM{Y,k}{X_k}, \xiM{Z,k}{X_k})$ defined in (\ref
{eqximy}),~(\ref{eqximz}) solves a {MWDP} with terminal condition $0$
and driver
$f_{\xi^*,k } (y,z):= f_k(X_k,\yM{k+1}{X_{k+1}}, \zM k{X_k}
)-f_k(X_k,y_{k+1}(X_{k+1}), z_k(X_k) )$.
Applying\break 
Proposition~\ref{propstability} with {$(Y_2,Z_2)\equiv0$ (so that
$L_{f_2} = 0$)} and using the Lipschitz continuity of $f_j(\cdot )$ yields
\begin{eqnarray*}
\bigl\llvert\xiM{Y,k} {X_k} \bigr\rrvert_2 & \le&
L_f \sum_{j =
k}^{N-1}
\frac{\bar\mathcal{E}
(Y,M,j+1) + \bar\mathcal{E}(Z,M,j) }{(T-t_j)^{1/2 - \sfrac
\thetaL2} } \Delta_j,
\\
\bigl\llvert\xiM{Z,k} {X_k} \bigr\rrvert_2 & \le&
C_{M}L_f \sum_{j =
k+1}^{N-1}
\frac{\bar
\mathcal{E}
(Y,M,j+1) + \bar\mathcal{E}(Z,M,j) }{(T-t_j)^{1/2 - \sfrac
\thetaL2}
\sqrt
{t_j - t_k} } \Delta_j.
\end{eqnarray*}
Next, introducing the notation $\Theta_j:= \mathcal{E}(Y,M,j+1) +
\mathcal{E}(Z,M,j)$
and applying Proposition~\ref{propeqyzerrdecoM}, it follows that
\begin{eqnarray*}
\bigl\llvert\xiM{Y,k} {X_k} \bigr\rrvert_2 & \le&\sqrt2
L_f \sum_{j = k}^{N-1}
\frac{
\Theta_j \Delta_j }{(T-t_j)^{1/2 - \sfrac\thetaL2} } + L_f \sum_{j
= k}^{N-1}
\frac{ ( \mathcal{E}^Y_{\mathrm{Dep.},j+1}
+ \mathcal{E}
^Z_{\mathrm{Dep.},j} ) \Delta_j
}{(T-t_j)^{1/2 - \sfrac\thetaL2} },
\\
\bigl\llvert\xiM{Z,k} {X_k} \bigr\rrvert_2 & \le&\sqrt2
C_{M}L_f \sum_{j = k +1}^{N-1}
\frac{
\Theta_j \Delta_j }{(T-t_j)^{1/2 - \sfrac\thetaL2} \sqrt{t_j -
t_k}}
\\
&&{} + C_{M}L_f \sum
_{j = k+1}^{N-1} \frac{ ( \mathcal{E}^Y_{\mathrm
{Dep.},j+1}
 +
\mathcal{E}^Z_{\mathrm{Dep.},j} )
\Delta_j
}{(T-t_j)^{1/2 - \sfrac\thetaL2} \sqrt{t_j - t_k} }.
\end{eqnarray*}
Substituting the above into (\ref{eqMCy3}) and (\ref{eqMCz1}),
and merging together the terms in $Z$, it follows that
%
\begin{eqnarray}
\mathcal{E}(Y,M,k)& \le&\mathcal{E}^Y_{\mathrm{App.},k} +\ObsBd{y,k}
\sqrt{\frac{
K_{Y,k}}{M_k}} + 2L_f \sum_{j = k}^{N-1}
\frac{ ( \mathcal{E}^Y_{\mathrm
{Dep.},j+1} + \mathcal{E}
^Z_{\mathrm{Dep.},j} ) \Delta_j
}{(T-t_j)^{1/2 - \sfrac\thetaL2} }\nonumber
\\
&&{} + 4 L_f \sum_{j = k}^{N-1}
\frac{ \Theta_j \Delta_j
}{(T-t_j)^{1/2 - \sfrac\thetaL2} }
\nonumber\\[-8pt]\label{eqMCy4} \\[-8pt]\nonumber
& \le&\mathcal{E}^Y_{\mathrm{App.},k} +\ObsBd{y,k} \sqrt{
\frac{
K_{Y,k}}{M_k}} + 2 \sum_{j = k}^{N-1}
\frac{ \mathcal{E}(j) \Delta_j
}{(T-t_j)^{1/2 - \sfrac\thetaL2} }
\\
&&{} + 4 L_f \sum_{j = k}^{N-1}
\frac{ \Theta_j \Delta_j
}{(T-t_j)^{1/2 - \sfrac\thetaL2} }, \nonumber
\\
\mathcal{E}(Z,M,k) & \le&\mathcal{E}^Z_{\mathrm{App.},k} +\ObsBd{z,k}
\sqrt{\frac
{K_{Z,k}} {
M_k} } + C_{M}\sum_{j = k+1}^{N-1}
\frac{ \mathcal{E}(j) \Delta_j
}{(T-t_j)^{1/2 - \sfrac\thetaL2} \sqrt{t_j - t_k} }
\nonumber\\[-8pt]\label{eqMCz2}\\[-8pt]\nonumber
&&{}+ \sqrt2 C_{M}L_f \sum
_{j = k +1}^{N-1} \frac{
\Theta_j \Delta_j }{(T-t_j)^{1/2 - \sfrac\thetaL2} \sqrt{t_j -
t_k}}.
\end{eqnarray}
\end{longlist}

\begin{longlist}
\item[\textit{Step} 4: \textit{final estimates.}]
Now, summing (\ref{eqMCz2}) and (\ref{eqMCy4}), one obtains an
estimate for $\Theta_k$:
\begin{eqnarray*}
\Theta_k & \le&\mathcal{E}(k) + (C_{M}+ 2\sqrt T) \sum
_{j =
k+1}^{N-1} \frac{
\mathcal{E}(j) \Delta_j }{(T-t_j)^{1/2 - \sfrac\thetaL2} \sqrt
{t_j - t_k}
}
\nonumber
\\
&&{}+ L_f( \sqrt2 C_{M}+ 4 \sqrt T) \sum
_{j = k +1}^{N-1} \frac{
\Theta_j \Delta_j }{(T-t_j)^{1/2 - \sfrac\thetaL2} \sqrt{t_j - t_k}}.
\end{eqnarray*}
Thus, using Lemmas~\ref{lemiterationgen1} and~\ref
{lemiterationgen2} with $\alpha= 0$, $\beta= \frac\thetaL2$, $C_u
= L_f( \sqrt2 C_{M}+ 4 \sqrt T)$, $w_k:=\mathcal{E}(k) + (C_{M}+
2\sqrt
T) \sum_{j = k+1}^{N-1} \frac{
\mathcal{E}(j) \Delta_j }{(T-t_j)^{1/2 - \sfrac\thetaL2} \sqrt
{t_j - t_k}
}$, we can control weighted sums involving $(\Theta_k)_k$ using
weighted sums of $(w_k)_k$, which is exactly what we need to complete
the upper bounds (\ref{eqMCy4})--(\ref{eqMCz2}) for $\mathcal
{E}(Y,M,k)$ and
$\mathcal{E}(Z,M,k)$. Namely, let $\gamma>0$:
\begin{eqnarray*}
&& \sum_{j = k+1}^{N-1} \frac{ w_j \Delta_j }{(T-t_j)^{1/2 - \sfrac
\thetaL2} (t_j - t_k)^{1-\gamma} }
\\
&&\quad  \le\sum
_{j = k+1}^{N-1} \frac{ \mathcal{E}(j) \Delta_j
}{(T-t_j)^{1/2 - \sfrac\thetaL2} (t_j - t_k)^{1-\gamma} }
\\
&&\qquad{}+(C_{M}+
2\sqrt T)\sum_{l=k+2}^{N-1} \frac{ \mathcal{E}(l) \Delta_l
}{(T-t_l)^{1/2 - \sfrac
\thetaL2} }
\sum_{j = k+1}^{l-1} \frac{ \Delta_j
}{(t_l-t_j)^{ 1 - \sfrac\thetaL2} (t_j - t_k)^{1-\gamma} }
\\
&&\quad  \le\bigl(1+\intB{\sfrac\thetaL2} {\gamma} T^{\thetaL/2}(C_{M}+ 2
\sqrt T)\bigr) \sum_{l=k+1}^{N-1}
\frac{ \mathcal{E}(l) \Delta_l}{
(T-t_l)^{1/2 - \sfrac\thetaL2} (t_l - t_k)^{1 - \gamma} },
\end{eqnarray*}
where we have applied Lemma~\ref{lemintegrals}. Thus,
\begin{eqnarray*}
&& \sum_{j = k+1}^{N-1} \frac{
\Theta_j \Delta_j
}{(T-t_j)^{1/2 - \sfrac\thetaL2} (t_j - t_k)^{1-\gamma} }
\\
&&\quad  \le
\cwg{\gamma}\bigl(1+\intB{\sfrac\thetaL2} {\gamma} T^{\sfrac\thetaL
2}(C_{M}+
2\sqrt T)\bigr) \sum_{l=k+1}^{N-1}
\frac{ \mathcal{E}(l) \Delta_l}{
(T-t_l)^{1/2 - \sfrac\thetaL2} (t_l - t_k)^{1 - \gamma} }
\end{eqnarray*}
and plugging the above inequality into (\ref{eqMCy4}) and (\ref
{eqMCz2}) yields (\ref{eqMCyerr}) and (\ref{eqMCzerr}).
\end{longlist}

\subsection{Complexity analysis}\label{sectioncomplexity}

As usual in empirical regression theory, appropriately tuning numerical
paramaters is crucial for finding the right trade-off between
statistical errors and estimation errors. This analysis allows to
express the error magnitude as a function of computational work
(complexity analysis).
We discuss { the complexity }in different cases according to the
regularity of the value functions $(y_i(\cdot), z_i(\cdot))$ and the
choice of the grid $\pi$. In order to have a fair comparison with other
numerical schemes, we revisit the setting of \cite{gobeturk14}, Section~4.4, which we partly recall for completeness,
and extend the analysis
{to include more general settings.}
\begin{itemize}
\item We perform an asymptotic complexity analysis as the number $N$ of
grid times goes to $+\infty$. We are concerned with time-dependent
bounds: thus in the following, the order convention, $\mathrm{O}(\cdot )$ or $\mathrm{o}(\cdot )$,
is uniform in $t_i$.
\item The grids under consideration are of the form $\pithetapi:=\{t_i
= T- T(1-\frac{i}{N})^{\sfrac{1}{\tpi}} \}$ for $\theta_\pi\in(0,1]$
(inspired by \cite{gobemakh10,geisgeisgobe12}). Observe that their
time-step $\Delta_i$ is not-increasing in $i$, hence they all satisfy
\HFc~with the same parameter $R_\pi=1$.

\item The magnitude of the final accuracy is denoted by $N^{-\thetaconv
}$ for some parameter $\thetaconv>0$. This is usually related to
time-discretization errors between the continuous-time BSDE and the
discrete-time one, $\thetaconv$ may range from $0^+$ (for non-smooth
data \cite{gobemakh10}, Theorem 1.1) to~1 (in the case of smooth data
\cite{gobelaba07}, Theorems 7 and 8).

\item The approximation\vspace*{1pt} spaces are given by local polynomials of degree
$n$ $(n\geq0)$ defined on hypercubes with edge length $\delta>0$,
covering the set $[-R,R]^d$ ($R>0$): we denote it by $\locpol{n,\delta
,R}$. The functions in $\locpol{n,\delta,R}$ take values in $\R$ for
the $y$-component and in $(\R^q)^\top$ for~$z$ {(using local
polynomials component-wise)}, but we omit this in the notation. The
best-approximation errors are easily controlled (using the Taylor formula):
%
\begin{equation}
\label{eqpolloc} \inf_{\varphi\in\locpol{n,\delta,R}}\bigl\llvert
\varphi
(X_i)-u(X_i)\bigr\rrvert_2\leq\llvert u
\rrvert_\infty\bigl(\P\bigl(\llvert X_i\rrvert
_\infty>R\bigr)\bigr)^{1/2}+ c_n \bigl\llvert
D^{n+1} u\bigr\rrvert_\infty\delta^{n+1}
\end{equation}
for any function $u$ that is bounded, $n+1$-times continuously
differentiable with bounded derivatives, and where the constant $c_n$
does not depend on $(R,u,\delta)$. The dimension of the vector space
$\locpol{n,\delta,R}$ is bounded by $\tilde c_n (2R/\delta)^{d}$ where
$\tilde c_n$ {is the number of polynomials on each hypercube (it
depends on $d$ and $n$)}.

A significant computational advantage of local polynomial basis is that
the cost of computing the regression coefficients associated to a
sample of size $M\geq\operatorname{dim}(\locpol{n,\delta,R})$ is $\mathrm{O}(M)$ flops.
{The cost of the regression in the $l$th hypercube is of order
$M^{(l)} \times{\tilde c_n^2}$ using SVD least squares minimization
\cite{goluvanl96}, Chapter~5, where $M^{(l)}$ is the number of
simulations that land in the hypercube.
Therefore, the total cost of the regressions at any time-point is of
order ${\tilde c_n^2} \sum_l M^{(l)} = {\tilde c_n^2} M = \mathrm{O}(M)$.
}

On the other hand, the cost of generating the clouds of simulations and
computing the simulated functionals $(\ObsM{Y,i}{X^{(i,m)}},\ObsM
{Z,i}{{H}^{(i,m)},{X}^{(i,m)}})_{i,m}$ is $\mathrm{O}(\sum_{i=0}^{N-1}NM_i)$,
which is clearly dominant in the computational cost $\mathcal{C}$ of
the MWLS
algorithm. To summarize, the computational cost is
\[
\mathcal{C}=\mathrm{O}\Biggl(\sum_{i=0}^{N-1}N
M_i\Biggr).
\]
{Another advantage of the local polynomial basis is that there is
substantial potential for parallel computing.}

\item To make the tail contributions (outside $[-R,R]^d$) small enough,
we assume that $X_i$ has exponential moments (uniformly in $i$), that
is, $\sup_{N\geq1}\sup_{0\leq i \leq N}\E(\mathrm{e}^{\lambda\llvert
X_i\rrvert_\infty
})<+\infty$ for some $\lambda>0$, so that the choice
$ R:= 2\thetaconv\lambda^{-1} \log(N+1)$ is sufficient to ensure
$(\P
(\llvert X_i\rrvert_\infty> R))^{1/2}= \mathrm{O}(N^{-\thetaconv})$.
\end{itemize}
To simplify the discussion, we assume $\thetaL=\thetaC=1$.

\subsubsection*{Smooth functions} Assume that $y_i(\cdot), z_i(\cdot)$ are,
respectively, of class $\mathcal{C}^{l+1}_b(\R^d,\R)$ and $\mathcal
{C}^{l}_b(\R
^d,(\R
^q)^\top)$ (bounded with bounded derivatives) for some $l \in\N
\setminus\{0\}$: this is similar to the discussion of \cite{gobeturk14}, Section~4.4. In fact, this is usually valid for the
continuous-time limit (a priori estimates on the semi-linear PDE, see
\cite{delaguat06,crisdela12}) provided that the data are smooth
enough. In particular, we may assume \HGt~with $\theta_\Phi=1$. This
leads to time-uniform bounds on the quantities $C_{y,i}, C_{z,i},
\ObsBd
{y,i}, \sqrt{T-t_i}\ObsBd{z,i}$.

Set
\[
\delta_{y,i}:=N^{-\afrac{\thetaconv}{l+1}},\qquad\delta_{z,i}:=N^{-\sfrac {\thetaconv}{l }},
\qquad M_i:=\bigl(\log(N+1)\bigr)^{d+1}
N^{\thetaconv
(2+\sfrac{ d }{l })},
\]
take $\mathcal{K}_{Y,i}:=\locpol{l,\delta_{y,i},R}$ and $\mathcal
{K}_{Z,i}:=\locpol
{l-1,\delta_{z,i},R}$.
From Proposition~\ref{propeqyzerrdecoM}, Theorem~\ref{thmMCerr}
and the inequality (\ref{eqpolloc}), it is easy to check that
\begin{eqnarray*}
\mathcal{E}^Y_{\mathrm{App.},i} &=&\mathrm{O}\bigl(N^{-\thetaconv}\bigr),
\qquad\mathcal{E}^Y_{\mathrm
{Dep.},i}=\mathrm{o}\bigl(N^{-\thetaconv
}\bigr),
\\
\ObsBd{y,i}\sqrt{\frac{ K_{Y,i}}{M_{i} }}&=&\mathrm{o} \bigl(N^{-\thetaconv
}/\sqrt{
\log(N+1)} \bigr),
\\
\mathcal{E}^Z_{\mathrm{App.},i}&=&\mathrm{O}\bigl(N^{-\thetaconv}\bigr),
\qquad\mathcal{E}^Z_{\mathrm{Dep.},i}=\mathrm{O}\bigl(N^{-\thetaconv}\bigr),
\\
\ObsBd{z,i}\sqrt{\frac{ K_{Z,i}}{M_{i} }}&=&(T-t_i)^{-\sfrac{ 1 }{2
}}\mathrm{O}
\bigl(N^{-\thetaconv}/\sqrt{\log(N+1)} \bigr).
\end{eqnarray*}
Consequently, using Lemma~\ref{lemintegrals}, we finally obtain
\begin{eqnarray*}
\bigl(\E\bigl[\bigl\llVert y_i - y^M_i
\bigr\rrVert_{i,M} ^2\bigr] \bigr)^{1/2}&=&\mathrm{O}
\bigl(N^{-\thetaconv
}\bigr),
\\
\bigl(\E\bigl[\bigl\llVert z_i -
z^M_i\bigr\rrVert_{i,M} ^2
\bigr] \bigr)^{1/2}&=&\mathrm{O}\bigl(N^{-\thetaconv
}\bigr) \biggl(1+
\frac{(T-t_i)^{-\sfrac{ 1 }{2 }}}{\sqrt{\log(N+1)}} \biggr).
\end{eqnarray*}
For\vspace*{1pt} any time-grid $\pi=\pithetapi$, we get $\sup_{0\leq i\leq N}\E
[\llVert y_i - y^M_i \rrVert_{i,M} ^2]+ \sum_{i=0}^{N-1}\Delta_i \E
[\llVert z_i - z^M_i\rrVert_{i,M} ^2]= \mathrm{O}(N^{-2\thetaconv})$.
The computational cost is
$\mathcal{C}=\mathrm{O} (\log(N+1)^{d+1} N^{\thetaconv(2+\sfrac{ d }{l
})+2})$.
Ignoring the logarithmic factors, we obtain a final accuracy in terms
of the computational cost:
\[
\mathcal{C}^{-\afrac{ 1 }{(2+\sfrac{d }{l })+ \sfrac{2 }{\thetaconv} }}.
\]
It should be compared with the rate $\mathcal{C}^{-\afrac{1}{(2+\sfrac
{d }{l })+
\sfrac{3 }{\thetaconv}}}$ which is valid for the Least Squares
Multi-step forward Dynamical Programming algorithm (LSMDP) \cite{gobeturk14}.
This shows an improvement on the rate, {although there is no change in
the dependence on dimension.}
The ratio $d/l$ is the usual balance between dimension and smoothness,
arising when approximating a multi-dimensional function.
{The controls of MWLS are stated in stronger norms than the controls of
LSMDP, and despite that, the estimates improve. The convergence
improvement is due to better MWDP-intrinsic estimates on $Z$, which
avoid the $1/\Delta_i$-factor of the LSMDP. This results in better
local error bounds, whence better global estimates. The reader can
easily check that this happens already in the simple case with null driver.}

\subsubsection*{H\"older terminal condition} We investigate the case of
non-smooth terminal condition, where nevertheless there is a smoothing
effect of the conditional expectation yielding smooth value functions
$(y_i(\cdot),z_i(\cdot))$. Namely, assume that $\Phi$ is bounded and
$\thetaP$-H\"older continuous (in particular with \HGt), and that,\vspace*{2pt} for
all~$i$, the function $y_i(\cdot)$ (resp., $z_i(\cdot)$) is
$(l+1)$-times (resp., $l$-times) continuously differentiable with
highest derivatives bounded by
%
\begin{eqnarray}
\bigl\llvert D^{l+1}_x y_i\bigr\rrvert
_\infty\le C(T-t_i)^{(\thetaP-
l)/2},\qquad\bigl\llvert
D^l_x z_i\bigr\rrvert_\infty\le
C(T-t_i)^{(\thetaP- (l+1))/2}. \label
{eqderivyz}
\end{eqnarray}
These qualitative assumptions are related to the works of \cite
{delaguat06,crisdela12}, who have determined similar estimates for
the gradients of quasi-linear PDEs under quite general conditions on
the driver, terminal condition and differential operator.
Their estimates cover the case $l=0$
(\cite{delaguat06}, Theorem~2.1) or $\thetaP=0$ {and $l \ge1$} (\cite{crisdela12}, Theorem 1.4), but the H\"older continuous setting {with
high order derivatives} is not investigated. We therefore extrapolate
these results in the assumptions (\ref{eqderivyz}) for the purposes
of this discussion.

In this setting, we have time-uniform bounds on the quantities $C_{y,i}$,
$(T-t_i)^{({ 1-\thetaP})/{2 }}C_{z,i}$, $\ObsBd{y,i}$, $\sqrt
{T-t_i}\ObsBd{z,i}$.
Set
\begin{eqnarray*}
\delta_{y,i}&:=&\sqrt{T-t_i}N^{-\afrac{\thetaconv}{l+1}},\qquad\delta
_{z,i}:=\sqrt{T-t_i}N^{-\sfrac{\thetaconv}{l }},
\\
M_i&:=&\bigl(\log(N+1)\bigr)^{d+1} N^{\thetaconv(2+\sfrac{ d }{l
})}(T-t_i)^{-d/2},
\end{eqnarray*}
take $\mathcal{K}_{Y,i}:=\locpol{l,\delta_{y,i},R}$ and $\mathcal
{K}_{Z,i}:=\locpol
{l-1,\delta_{z,i},R}$.
Similarly to before, using in particular (\ref{eqpolloc}), we
eventually obtain
\begin{eqnarray*}
\mathcal{E}^Y_{\mathrm{App.},i} &=&\mathrm{O}\bigl(N^{-\thetaconv}\bigr),
\qquad\mathcal{E}^Y_{\mathrm
{Dep.},i}=\mathrm{o}\bigl(N^{-\thetaconv
}\bigr),
\\
\ObsBd{y,i}\sqrt{\frac{ K_{Y,i}}{M_{i} }}&=&\mathrm{o} \bigl(N^{-\thetaconv
}/\sqrt{
\log(N+1)} \bigr),
\\
\mathcal{E}^Z_{\mathrm{App.},i}&=&(T-t_i)^{\vfrac{\thetaP- 1 }{2
}}\mathrm{O}
\bigl(N^{-\thetaconv}\bigr), \qquad\mathcal{E}^Z_{\mathrm
{Dep.},i}=(T-t_i)^{\vfrac{\thetaP- 1 }{2
}}\mathrm{O}
\bigl(N^{-\thetaconv}\bigr),
\\
 \ObsBd{z,i}\sqrt{\frac{ K_{Z,i}}{M_{i} }}&=&(T-t_i)^{-\sfrac{ 1 }{2
}}\mathrm{O}
\bigl(N^{-\thetaconv}/\sqrt{\log(N+1)} \bigr).
\end{eqnarray*}
Consequently, using Lemma~\ref{lemintegrals}, we finally obtain
\begin{eqnarray*}
 \bigl(\E\bigl[\bigl\llVert y_i - y^M_i
\bigr\rrVert_{i,M} ^2\bigr] \bigr)^{1/2}&=&\mathrm{O}
\bigl(N^{-\thetaconv
}\bigr),
\\
 \bigl(\E\bigl[\bigl\llVert z_i - z^M_i
\bigr\rrVert_{i,M} ^2\bigr] \bigr)^{1/2}&=&\mathrm{O}
\bigl(N^{-\thetaconv
}\bigr) \biggl((T-t_i)^{\vfrac{\thetaP- 1 }{2 }}+
\frac{(T-t_i)^{-\sfrac{ 1 }{2
}}}{\sqrt
{\log(N+1)}} \biggr).
\end{eqnarray*}
The computation cost is given by (under the assumption $\pi=\pithetapi$)
\[
\mathcal{C}=\mathrm{O} \Biggl(\sum_{i=0}^{N-1}N
M_i \Biggr)=\mathrm{O} \bigl(\bigl(\log(N+1)\bigr)^{d+1}
N^{1+\thetaconv(2+\sfrac{ d }{l })} \bigr) \sum_{i=0}^{N-1}
\biggl(1-\frac{ i }{ N}\biggr)^{-\afrac{ d }{2 \tpi}}.
\]
Up to possibly a $\log(N)$-factor, the last sum is $\mathrm{O}(N^{\afrac{ d }{2
\tpi}\lor1})$. Ignoring the logarithmic factors, we obtain
$\mathcal{C}=\mathrm{O}(N^{1+\afrac{ d }{2 \tpi}\lor1+\thetaconv(2+\sfrac{d
}{l })})$.
Equivalently, as a function of the computational cost, the convergence
rate of the final accuracy equals
\[
\mathcal{C}^{-\sfrac{ 1 }{ ((2+\sfrac{d }{l })+\sfrac{1 }{\thetaconv
}(1+\sfrac{ d}{(2 \tpi)}\lor1)})}.
\]
Following \cite{gobemakh10} (under suitable assumptions), two
time-grid choices are possible for solving the same BSDE.
\begin{itemize}
\item The uniform grid $\pi=\pit{1}$ gives $\thetaconv=\thetaP/2$ (at
least). The convergence order becomes
{$( 2+\frac{ d }{l }+\frac{ 2 }{\thetaP}(1+ \frac{d }{2 } \lor1 ) )^{-1}$.}

\item The grid $\pi=\pit{\theta}$ (for $\theta<\thetaP$) gives
$\thetaconv=1/2$. Taking $\theta\uparrow\thetaP$, the convergence
order is
{$(2+\frac{ d }{l }+\frac{2}{\thetaP}(\thetaP+ \frac{d }{2 } \lor
\thetaP
))^{-1}$.}
\end{itemize}
The grid $\pit{\theta}$ exhibits a better convergence rate compared to
the uniform grid. This corroborates the interest in time grids that are
well adapted to the regularity of the data.
These features will be investigated in subsequent more experimental works.

\begin{appendix}
\section*{Appendix}\label{appendix}

\setcounter{subsection}{0}
\setcounter{equation}{0}

\subsection{Proof of Lemmas \texorpdfstring{\protect\ref{lemintegrals}}{2.1}, \texorpdfstring{\protect\ref{lemiterationgen1}}{2.2}
and \texorpdfstring{\protect\ref{lemiterationgen2}}{2.3}}\label{appendixprooflemintegrals}

\subsubsection{Proof of Lemma \texorpdfstring{\protect\ref{lemintegrals}}{2.1}}
The first inequality, for $\alpha\leq1$, follows by bounding the sum
by $\int_{t_i}^{t_k} (t_k-t)^{\alpha-1} \,\mathrm{d}t$, whence $\intB{\alpha}{1}
=1/\alpha$. The case $\alpha>1$ is obvious with $\intB{\alpha}{1}=1$.
For the second inequality, there are two main cases:
\begin{longlist}
\item[$\rhd$] If $\alpha\geq1$ and $\beta\geq1$, the advertised inequality
is obvious with $B_{ \alpha,\beta}=1$.

\item[$\rhd$] Now, assume the complementary case, that is, $\alpha< 1$
and/or $\beta< 1$, and first
consider the case $t_{i} = 0$ and $t_{k} = 1$.
We set $\varphi(s) = (1-s)^{\alpha- 1}s^{\beta-1}$ and we
use the integral $\int_{0}^1\varphi(s) \,\mathrm{d}s$
(equivalent to the usual beta function with parameters $(\alpha,\beta)$)
to bound the sum. {A simple but useful property {(due to $\alpha< 1$
and/or $\beta< 1$)} is that $\varphi$ is either monotone or {has a
unique minimum on $(0,1)$}, whence
%
\begin{eqnarray*}
(1 - t_{j})^{\alpha-1}t_{j}^{\beta-1}
\Delta_{j} \le R_\pi\int_{t_{j-1}}^{t_{j}}
\varphi(s) \,\mathrm{d}s+\int_{t_j}^{t_{j+1}}\varphi(s) \,\mathrm{d}s.
\end{eqnarray*}
Summing up over $j$ and defining $B_{\alpha,\beta} = (1+R_{\pi})
\int_{0}^1 \varphi(s) \,\mathrm{d}s$ concludes the proof for the simple case.}
For general $t_{i}$ and $t_{k}$ one can use the bounds on the simple
case by rearranging the $j$-sum which is equal to
\begin{eqnarray*}
(t_{k} - t_{i})^{\alpha+ \beta-1} \sum
_{j=i+1}^{k-1}\biggl( 1 - \frac{t_{j}-t_{i}}{t_{k} - t_{i}}
\biggr)^{\alpha-1} \biggl(\frac{t_{j} - t_{i}}{t_{k} - t_{i}}\biggr
)^{\beta-1}
\frac{\Delta_{j}}{t_k - t_i} \le B_{\alpha,\beta} (t_{k} -
t_{i})^{\alpha+ \beta-1}.
\end{eqnarray*}
\end{longlist}
%

\subsubsection{Proof of Lemma \texorpdfstring{\protect\ref{lemiterationgen1}}{2.2}}
If $\alpha\geq\frac{1 }{2 }$, the result trivially holds with $\cw=1$
and $\cuh=C_{u}{T^{\alpha-\sfrac{ 1}{2 }}}$.

Now, assume $\alpha<\frac{1 }2 $: if (\ref
{eqiterationfeed}) holds, of course we also have
%
\begin{equation}\label{eqiterationgenbis}
\hspace*{-10pt}u_j \le w_j + \sum_{l=j+1}^{N-1}
\frac{w_l \Delta_l}{(T-t_l)^{{{1}/2 - \beta}}(t_l - t_j)^{{{1}/2-\alpha}}}
+C_{u}\sum_{l=j+1}^{N-1}
\frac{u_l \Delta_l}{(T-t_l)^{{{1}/2 -
\beta}}(t_l - t_j)^{{{1}/2-\alpha}}}.
\end{equation}
By substituting (\ref{eqiterationgenbis}) into the last sum, and
using Lemma~\ref{lemintegrals} we observe
\begin{eqnarray*}
&&\sum_{l=j+1}^{N-1} \frac{u_l \Delta_l }{(T-t_l)^{\sfrac{1} 2- \beta}(t_l
- t_j)^{\sfrac{1} 2- \alpha}}
\\
&&\quad \le\sum_{l=j+1}^{N-1} \frac{w_l\Delta_l}{(T-t_l)^{\sfrac{1 }{2 }-
\beta}(t_l - t_j)^{\sfrac{1} 2- \alpha}}
\\
&&\qquad{} +
\sum_{l = j+1}^{N-1} \frac{\sum_{r = l+1}^{N-1} {w_r\Delta
_r}\Delta_l
/({(T-t_r)^{1/2-\beta}(t_r-t_l)^{1/2 -\alpha}}) }{
(T-t_l)^{\sfrac{1} 2- \beta}(t_l - t_j)^{\sfrac{1} 2- \alpha}}
\\
&&\qquad{}+ C_{u}\sum_{l = j+1}^{N-1}
\frac{\sum_{r = l+1}^{N-1}
{u_r\Delta_r}
\Delta_l/
({(T-t_r)^{\sfrac{1} 2- \beta}(t_r - t_l)^{\sfrac{1} 2-
\alpha
}}) }{
(T-t_l)^{\sfrac{1} 2- \beta}(t_l - t_j)^{\sfrac{1} 2- \alpha}}
\\
&&\quad \le\sum_{l=j+1}^{N-1} \frac{w_l\Delta_l}{(T-t_l)^{\sfrac{1 }{2 }-
\beta}(t_l - t_j)^{\sfrac{1} 2- \alpha}}
\\
&&\qquad{}+
B_{\alpha+\beta,\sfrac{1}2+\alpha} \sum_{r = j+2}^{N-1}
\frac
{w_r\Delta
_r}{(T-t_r)^{1/2-\beta}(t_r-t_j)^{1/2 -2\alpha-\beta}}
\\
&&\qquad{}+ C_{u}B_{\alpha+\beta,\sfrac{1}2+\alpha} \sum_{r =
j+2}^{N-1}
\frac{u_r\Delta_r}{(T-t_r)^{1/2-\beta
}(t_r-t_j)^{1/2 -2\alpha-\beta}}.
\end{eqnarray*}
Substituting into (\ref{eqiterationgenbis}), we observe that we have
an equation of similar form to (\ref{eqiterationgenbis}),
except that, in the sum involving $u$, $\alpha\mapsto2\alpha+ \beta$
and $C_u \mapsto C_u^2 B_{\alpha+\beta,\sfrac{ 1 }{2 }+\alpha}$,
and, in the sum involving $w$, $w \mapsto{(1+C_u (1+ T^{\alpha+\beta
}B_{\alpha+\beta,\sfrac{ 1 }{2 }+\alpha} ))}w$.

After $\kappa$ iterations of the previous step, we obtain $\alpha
\mapsto2^{\kappa}(\alpha+\beta) - \beta=: \alpha_\kappa$.
Hence, for $\kappa$ sufficiently large so that $\alpha_\kappa\ge
\frac{1}{2}$,
that is, $\kappa\ge\log_2 (\frac{\sfrac{1}2 + \beta}{\alpha+
\beta }) $,
we obtain the bound advertised in the lemma statement.

\subsubsection{Proof of Lemma \texorpdfstring{\protect\ref{lemiterationgen2}}{2.3}}
{W.l.o.g. we can assume that $\cw=1$ in (\ref{eqiterationfeed2});
if it is not, one can redefine $w$ as $\cw w$.}
We first prove the case $\gamma=1$.
Define
%
\begin{equation}
\zeta_{s}:= 2\cuh\int_{0}^{s}
\frac{\mathrm{d}r}{(T - r)^{1/2 -\beta}} \le\frac{ 2 }{1+2\beta} 2\cuh
{T^{(1+2\beta)/2}},
\label{eqbeta}
\end{equation}
and write $\zeta_{j } = \zeta_{t_{j}}$ for brevity.
{We first multiply (\ref{eqiterationfeed2}) by ${\mathrm{e}^\zeta_j \Delta_j
\over(T-t_j)^{1/2 - \beta} }$, then sum the outcome equation over
$j \in\{i+1, \ldots, N-1\}$, and finally switch the order of
summation on the right-hand side as follows:}
\begin{eqnarray*}
&& \sum_{j=i+1}^{N-1} \frac{u_j \mathrm{e}^{\zeta_{j}}\Delta
_{j}}{(T-t_{j})^{1/2 - \beta} }
\\
&&\quad  \le
\sum_{j=i+1}^{N-1} \frac{w_j \mathrm{e}^{\zeta_{j}}\Delta
_{j}}{(T-t_{j})^{1/2 - \beta} }
\\
&&\qquad{} +\sum
_{j=i+1}^{N-1} \frac{\sum_{l=j+1}^{N-1} {w_l \Delta
_l}\mathrm{e}^{\zeta_{j}}\Delta_{j}/
({(T-t_l)^{{{1}/2 - \beta}}(t_l - t_j)^{{{1}/2-\alpha}}})
}{(T-t_{j})^{1/2 - \beta} }
\\
&&\qquad{} +\cuh\sum_{j=i+1}^{N-1} \frac{\sum_{l=j+1}^{N-1}
{u_l \Delta
_l}
\mathrm{e}^{\zeta_{j}}\Delta_{j}/
({(T-t_l)^{{{1}/2 - \beta}}})
}{(T-t_{j})^{1/2 - \beta} }
\\
&&\quad  \le \mathrm{e}^{\zeta_T}\sum_{j=i+1}^{N-1}
\frac{w_j \Delta
_{j}}{(T-t_{j})^{1/2 - \beta} }+\mathrm{e}^{\zeta_T}\intB{\alpha+\beta}
{1}\sum
_{l=i+2}^{N-1} \frac{w_l \Delta_l}{(T-t_l)^{{{1}/2 - \beta}}(t_l -
t_i)^{-\alpha
-\beta}}
\\
&&\qquad{}+\cuh\sum_{l=i+2}^{N-1} \frac{u_l\Delta_l}{(T-t_l)^{\sfrac{1} 2 -
\beta
}}
\sum_{j=i+1}^{l-1} \frac{ \mathrm{e}^{\zeta_{j}}\Delta_{j}}{(T-t_{j})^{1/2
- \beta} }
\\
&&\quad  \le \mathrm{e}^{\zeta_T}\bigl(1+\intB{\alpha+\beta} {1}T^{\alpha+\beta}\bigr)
\sum_{l=i+1}^{N-1} \frac{w_l \Delta_l}{(T-t_l)^{{{1}/2 - \beta}}} +
\frac{1}2 \sum_{l=i+1}^{N-1}
\frac{u_l \mathrm{e}^{\zeta_{l}}\Delta
_{l}}{(T-t_{l})^{1/2 - \beta} },
\end{eqnarray*}
where we have used (because $\zeta$ is {non-decreasing} and $\beta
\leq
\frac{1}2$)
\[
{\cuh}\sum_{j=i+1}^{l-1}\frac{ \mathrm{e}^{\zeta_{j}}\Delta
_{j}}{(T-t_{j})^{1/2 - \beta} }
\le\int_{{{t_{i+1}}}}^{t_l}\frac{{\cuh}\mathrm{e}^{\zeta_s }}{(T-s)^{1/2 -
\beta
}}\,\mathrm{d}s\leq{
\frac{\mathrm{e}^{\zeta_l}}{2}}. 
\]
{By subtracting} the term with factor $\frac{1}2$, the result for
$\gamma
=1$ follows.
Moreover, plugging the result into (\ref{eqiterationfeed2}), and
returning to general $\cw$, gives
%
\begin{eqnarray}\label{equpperboundu}
u_j &\le&\cuv w_j + \cuv\sum
_{l=j+1}^{N-1} \frac{w_l \Delta
_l}{(T-t_l)^{{{1}/2 - \beta}}(t_l - t_j)^{{{1}/2-\alpha}}}
\nonumber\\[-8pt]\\[-8pt]\nonumber
&&{} +\cuv\sum
_{l=j+1}^{N-1} \frac{w_l \Delta_l}{(T-t_l)^{{{1}/2 -
\beta
}}}
\end{eqnarray}
for a constant $\cuv:= 2 \cw \mathrm{e}^{\zeta_T} (1 + \intB{\alpha+ \beta}1
T^{\alpha+\beta})$. Now for the general case $\gamma>0$, observe that,
for any $\delta\geq0$,
one obtains by change of the order of summation that
%
\begin{eqnarray}\label{eqineqappendix}
&& \sum_{j=i+1}^{N-1}
\frac{\sum_{l=j+1}^{N-1}
{w_l \Delta_l}
\Delta_{j}/
({(T-t_l)^{{{1}/2 - \beta}}(t_l - t_j)^{{{1}/2-\delta}}}) }{(T-t_{j})^{1/2 - \beta
}(t_j-t_i)^{1-\gamma} }
\nonumber\\[-8pt]\\[-8pt]\nonumber
&&\qquad \leq B_{\beta+\delta,\gamma} \sum_{l=i+2}^{N-1}
\frac{w_l\Delta
_l}{(T-t_l)^{\sfrac{1} 2 - \beta}(t_l-t_i)^{1-\beta-\delta-\gamma} }.
\end{eqnarray}
Thus, (\ref{equpperboundu}) yields
\begin{eqnarray*}
&& \sum_{j=i+1}^{N-1} \frac{u_j \Delta_{j}}{(T-t_{j})^{1/2 - \beta
}(t_j -t_i)^{1-\gamma} }
\\
&&\quad \leq
\cuv\sum_{j=i+1}^{N-1} \frac{w_j \Delta_j}{(T-t_j)^{{{1}/2 -
\beta}}(t_j - t_i)^{{1-\gamma}}}
\\
&&\qquad{} +\cuv B_{\beta+\alpha,\gamma}\sum_{l=i+2}^{N-1}
\frac
{w_l\Delta_l}{(T-t_l)^{\sfrac{1} 2 - \beta}(t_l-t_i)^{1-\beta-\alpha
-\gamma} }
\\
&&\qquad{} +\cuv B_{\beta+\sfrac{1}2,\gamma}\sum_{l=i+2}^{N-1}
\frac
{w_l\Delta_l}{(T-t_l)^{\sfrac{1} 2 - \beta}(t_l-t_i)^{1/2-\beta
-\gamma} }
\\
&&\quad \leq\cuv\bigl(1+B_{\beta+\alpha,\gamma}T^{\alpha+\beta}+
B_{\beta+\sfrac{1}2,\gamma}T^{1/2+\beta}\bigr)\sum_{j=i+1}^{N-1}
\frac
{w_j \Delta_j}{(T-t_j)^{{{1}/2 - \beta}}(t_j - t_i)^{{1-\gamma}}}.
\end{eqnarray*}
%

\subsection{Proof of Lemma \texorpdfstring{\protect\ref{lemobsbounds}}{3.7}}\label{appendixprooflemmalemobsbounds}
Using the bounds $\Cyi$ and $\Czi$ on $\yM i\cdot$ and $\zM i\cdot$,
respectively, one applies the local Lipschitz continuity and
boundedness properties of $f_j$ given in \HF~to obtain the bound
%
\begin{equation}
\bigl\llvert f_j \bigl(x_j, \yM{j+1} {x_{j+1}}, \zM j{x_j} \bigr) \bigr\rrvert\le
\frac{L_f (\Cyjp
+ \Czj) }{(T-t_j)^{1/2 {- \sfrac\thetaL2}} } +
\frac
{C_f}{(T-t_j)^{1-\thetaC}}. \label{eqMCbdf}
\end{equation}
Substituting this into the definition $\ObsM{Y,i} {\vecx{x}^{(i)}} $
{(see (\ref{eqPsiMy}))}, it follows from \HF~that
\begin{eqnarray*}
\bigl\llvert\ObsMa{Y,i} {\vecx x^{(i)}}\bigr\rrvert\le
C_\xi+ \sum_{j = i}^{N-1} \biggl(
\frac{L_f (\Cyjp+ \Czj) }{(T-t_j)^{1/2 - \sfrac\thetaL2} } + \frac
{C_f}{(T-t_j)^{1-\thetaC}} \biggr)\Delta_j.
\end{eqnarray*}
Substituting the value of $\Cyj$ and $\Czj$ given in equations (\ref{eqasy}) and (\ref{eqasz}), respectively, using the crude bound
$\llvert\xi-\E_i\xi\rrvert_{2,i}\leq C_\xi$ and Lemma
\ref{lemintegrals}, we
obtain the bound $\ObsBd{y,i}$, with the form~(\ref{eqboundObsBd}).\vspace*{1pt}

To obtain the bound $\ObsBd{z,i}$, apply first the triangle inequality
on the conditional standard deviation of $\ObsM{l,Z,i}{\MW{i,m}{},
\X{i,m}{}}$; second use the bound (\ref{eqMCbdf}) on the driver,
and the bound \HH~to obtain
\begin{eqnarray*}
&& \sqrt{\var\bigl[\ObsMa{l,Z,i} {\MW{i,m} {}, \X{i,m} {}} \mid\F{M}i \bigr]}
\\
&&\quad \leq\frac{{C_\xi} C_{M}}{\sqrt{T - t_i}} + \sum_{j = i+1}^{N-1}
\biggl(\frac{L_f (\Cyjp+ \Czj)
}{(T-t_j)^{1/2 {- \sfrac\thetaL2}} } + \frac
{C_f}{(T-t_j)^{1-\thetaC
}} \biggr) \frac{ C_{M}}{\sqrt{t_j - t_i}}
\Delta_j.
\end{eqnarray*}
Then, the computation of $\ObsBd{z,i}$ follows again from equations
(\ref{eqasy}) and (\ref{eqasz}), and Lemma~\ref{lemintegrals}.
The form (\ref{eqboundObsBd}) is also derived. We skip details.
\end{appendix}

\section*{Acknowledgements}
This research benefited from the support of the Chair Finance and
Sustainable Development,
under the aegis of Louis Bachelier Finance and Sustainable Growth
laboratory, a joint initiative with Ecole polytechnique.
Emmanuel Gobet's
research is part of the Chair Financial Risks of the Risk Foundation
and of the FiME Laboratory.
Plamen Turkedjiev
acknowledges support form German Science Fondation DFG via Berlin
Mathematical School and Matheon.


%

\printhistory
\end{document}